\documentclass[a4paper,11pt]{article}

\usepackage{authblk}
\usepackage{amsmath}
\usepackage{amsfonts}
\usepackage{amssymb}
\usepackage{latexsym}
\usepackage{indentfirst}
\usepackage[top=1in, bottom=1in, left=1in, right=1in]{geometry}
\usepackage{algorithm}
\usepackage{algpseudocode}
\usepackage{graphicx,psfrag,epsf}
\usepackage{xcolor}
\usepackage[font=small]{caption}
\usepackage{subcaption} 
\usepackage{float}    
\usepackage{enumerate}
\usepackage{enumitem}
\usepackage{natbib}
\usepackage{url}
\usepackage{bm}
\usepackage{longtable}
\usepackage{multirow}
\usepackage{booktabs, adjustbox, multirow, makecell}
\usepackage{doi}

\let\footnote=\endnote

\newtheorem{theorem}{Theorem}

\newtheorem{remark}{Remark}
\newtheorem{corollary}{Corollary}
\newtheorem{proposition}{Proposition}
\newtheorem{assumption}{Assumption}

\allowdisplaybreaks

 \DeclareMathOperator*{\Var}{Var}
\DeclareMathOperator*{\EE}{E}

\DeclareMathOperator*{\dimm}{dim}

\def\om{\mathbf{\om}}
\def\si{\mathbf{\si}}

\newcommand{\ud}{\mathrm{d}}

\def\A{\mathbf{A}}

\def\F{\mathbf{F}}

\def\I{\mathbf{I}}

\def\K{\mathbf{K}}

\def\W{\mathbf{W}}
\def\X{\mathbf{X}}
\def\Y{\mathbf{Y}}
\def\Z{\mathbf{Z}}

\def\RR{\mathbb{R}}
\def\NN{\mathbb{N}}

\def\EE{\mathbb{E}}
\def\PP{\mathbb{P}}

\def\bfa{\mathbf{a}}

\def\bfk{\mathbf{k}}

\def\bfx{\mathbf{x}}
\def\bfy{\mathbf{y}}
\def\bfz{\mathbf{z}}

\begin{document}

\title{\bf Scalable Bayesian Quadrature via Operator-Adapted Features}

\author[1]{Huanyan Zhu}
\author[1]{Cheng Li}
\affil[1]{Department of Statistics and Data Science, National University of Singapore}
\date{}
\maketitle

\begin{abstract}
Kernel methods face two recurrent bottlenecks: the cubic cost of inverting an $n\times n$ kernel matrix, as in kernel ridge regression, and the evaluation of kernel integrals, as required in kernel quadrature. Existing low-rank kernel and Gaussian process approximations primarily address the former bottleneck, but do not control the downstream quadrature error when both the kernel matrix and kernel mean are approximated, particularly when the integrand lies outside the kernel reproducing kernel Hilbert space. We address both challenges through an operator-adapted finite-feature approximation. Specifically, we construct a finite-dimensional working kernel from the eigenvalues and eigenfunctions of a second-order differential operator. We establish convergence guarantees for the resulting finite-feature estimators in kernel quadrature, with extension to other related kernel methods. Our results show that, with $m \lesssim n$ basis functions, the truncated estimator attains the corresponding benchmark statistical rates. The theory covers fixed/random designs, noisy/noiseless observations, oversmoothing/undersmoothing regimes, and Mat\'ern/Gaussian-type weights. We provide explicit error bounds that quantify the roles of the spectral cutoff, regularization parameter, lengthscale, sampling geometry, and operator boundary condition. We verify the convergence results in numerical experiments.
\end{abstract}


\section{Introduction} \label{sec:intro}

Kernel methods provide a flexible and mathematically tractable framework for nonlinear estimation, prediction, and numerical integration. By representing functions through a positive definite kernel and the associated reproducing kernel Hilbert space (RKHS), they convert nonlinear problems on the input space into linear problems in a feature space, while retaining a rich approximation theory. This perspective underlies a wide range of methods in statistics, machine learning, spatial modeling, and probabilistic numerics, including Gaussian process regression, kernel ridge regression, support vector machines, and kernel quadrature. 

Among these, kernel quadrature is a statistically powerful tool that estimates the linear functional
\begin{align*}
    I(f) = \int_\Omega f(\bfx) \ud \PP(\bfx),
\end{align*}
for a function $f$ defined on the domain $\Omega \subset \RR^d$ equipped with a probability measure $\PP$ by exploiting smoothness through a positive-definite kernel $k(\cdot,\cdot)$. 
Under its Gaussian process (GP) interpretation, usually referred to as Bayesian quadrature (BQ) \cite{OHagan:1991BHquadrature, Briol:2019BQreview}, the same construction also provides a posterior distribution for $I(f)$. The estimator is the posterior mean taking the form
\begin{align} \label{eq:bq_predictor}
    \widehat I(f)
    =
    \bm\mu_n^\top (\K_n+\eta \I_n)^{-1} \Y_n,
\end{align}
where $\bm\mu_n = \left[ \int_\Omega k(\bfx_1, \bfx) ~\ud\PP(\bfx),\ldots, \int_\Omega k(\bfx_n,\bfx) ~\ud\PP(\bfx) \right]^\top$, $\K_n=\{k(\bfx_i,\bfx_j)\}_{1\le i,j\le n}$ is the $n\times n$ kernel matrix, $\I_n$ is the $n\times n$ identity matrix, $\Y_n$ is the vector of observations, and $\eta > 0$ is a regularization parameter. 
Exact kernel quadrature can attain high-order convergence rates for smooth integrands, including settings in which the target smoothness differs from that of the kernel \cite{Bach:2017, Kanagawa:2020}.
For Mat\'ern or Sobolev-equivalent kernels, Kanagawa et al. \cite{Kanagawa:2020} derived convergence rates for deterministic Bayesian quadrature under quasi-uniform designs, including misspecified settings in which the integrand is rougher than the imposed kernel RKHS. 
Wynne et al. \cite{wynne:2021} subsequently derived Bayesian-quadrature consequences for observations with correctly specified Gaussian noise. 

These statistical advantages, however, rely on the ability to compute both the quadrature weights and the kernel integral accurately and efficiently. Therefore, the implementation faces two distinct but intertwined bottlenecks: solving the $n\times n$ regularized kernel system generally requires cubic time cost, while evaluating $\bm\mu_n$ requires the kernel mean embedding to be analytically available or numerically tractable.

A natural approach is to replace the kernel by a finite-feature expansion, as in reduced-rank Gaussian process and spectral Bayesian quadrature methods \cite{RiutortMayol:2023, Solin:2020, Warren:2022, Warren:2024FFBQ}. Such an expansion reduces kernel matrix operations to an $m$-dimensional feature space and replaces the kernel mean by finitely many feature integrals. Nevertheless, a kernel-level approximation error analysis does not by itself yield a similar version for the resulting quadrature estimator, since the standard worst-case-error arguments for RKHS functions cannot apply when the integrand lies outside the RKHS. The problem becomes more delicate for feature systems constructed independently of the integration measure, and for Gaussian kernels, whose rapid spectral decay produces an exponential interaction among the feature cutoff, lengthscale, and regularization level. These difficulties lead to three research questions: 
\begin{itemize}
    \item In the misspecified case where $\nu>0$ denotes
    the Mat\'ern smoothness parameter, $H^s(\Omega)$ denotes the Sobolev space of order $s$ on $\Omega$, and $f \in H^s(\Omega)$ with $s < \nu + d/2$, what convergence rate can be established for the quadrature risk? Can this rate attain the benchmark rate?
    \item How do the number of features $m$ and the regularization parameter $\eta$ jointly affect the quadrature error? Can one obtain sharp statistical guarantees while reducing the computational complexity to $nm^2 + m^3$, $m \lesssim n$?
    \item What structural properties should a feature system satisfy to achieve an optimal convergence rate? Beyond Fourier features, which alternative feature families permit tractable feature integration, stable finite-dimensional inference, sharp quadrature guarantees, and calibrated uncertainty quantification?
\end{itemize}

\subsection{Main results and contributions}
We answer these questions by developing an operator-adapted finite-feature framework whose central analytical object is the quadrature estimator \eqref{eq:bq_predictor}, rather than the approximation error of the kernel or Gaussian process posterior mean. Let $\mathcal T_L$ be a proper second-order differential operator with compact resolvent on an auxiliary domain $D_L\Supset\Omega$, i.e., the closure $\overline \Omega$ is a compact subset of $D_L$ whose radius is $L$, and let $(\lambda_j,\phi_j)_{j\geq1}$ denote the eigenpairs of $\mathcal T_L$. We construct the finite-feature kernel
\begin{align*}
    \widehat k_m(\bfx,\bfx')
    :=
    \sum_{j=1}^{m}
    w(\lambda_j)
    \phi_j(\bfx)\phi_j(\bfx'),
    \qquad
    \bfx,\bfx'\in\Omega,
\end{align*}
where the spectral weights $w(\lambda_j)$ are chosen to reproduce the Mat\'ern or Gaussian type spectral decay. Unlike a Mercer expansion, the basis is constructed independently of $\PP$. Consequently, neither a truncated-kernel bound nor a posterior-mean approximation bound directly controls the resulting quadrature rule. Our analysis instead controls the regularized finite-feature weights and the quadrature error directly, separating the contributions of spectral truncation, regularization, observation noise, sampling geometry, kernel lengthscale, and the auxiliary domain $D_L$. At the computational level, the construction replaces the original kernel mean by the finite collection of feature integrals $\int_\Omega\phi_j(\bfx) \ud\PP(\bfx)$ and moves the matrix calculation from the $n$-dimensional observation space to the $m$-dimensional feature space.

Our first contribution is an end-to-end convergence theory for the resulting quadrature estimator over Sobolev classes $H^s(\Omega)$, with $s>d/2$, without requiring the integrand to belong to the kernel RKHS. For Mat\'ern spectral weights with $\tau=\nu+d/2$, the squared quadrature error under a quasi-uniform fixed design is of order $n^{-2\min\{s,\tau\}/d}$ in the noiseless case; under an i.i.d.\ random design, the same rate holds up to logarithmic factors with high probability. When $s\leq\tau$, this recovers the optimal deterministic Sobolev quadrature rate, while when $s>\tau$, it recovers the $\tau$-smooth benchmark established for exact Mat\'ern kernel quadrature \cite{Kanagawa:2020}. 
For noisy observations, both fixed and random designs attain the parametric squared-error order 
$n^{-1}$. 
Equivalently, the root mean squared error is $n^{-1/2}$. 
This result is sharper than $n^{-\frac{s}{2s+d}}$, which is the error rate for full-kernel quadrature \cite{wynne:2021}.

For Gaussian spectral weights, the noiseless squared-error rate attains the optimal Sobolev polynomial rate $n^{-2s/d}$, up to logarithmic factors. Our theory covers fixed and random designs, noisy and noiseless observations, Mat\'ern and Gaussian kernels, and both oversmoothing and undersmoothing regimes. The corresponding convergence rates and one set of sufficient feature dimensions are summarized in Table~\ref{tab:bq-rates}. 
For the bounded density class, the feature-budget conditions are not only sufficient but also necessary.

\begin{table*}[t]
\centering
\caption{Squared error rates for kernel quadrature with Laplacian/Legendre basis.
Here $\tau=\nu+d/2$. The error is defined as $\sup_{\|f\|_{H^s(\Omega)} \leq 1}\EE_\epsilon\Big[\big(I(f) - \widehat I_m(f) \big)^2 \mid \X_n \Big]$, where $\widehat I_m(f)$ is our finite-feature approximation with $m$ features and $\X_n$ is the design with $n$ points.
For random design, the error rate should be treated as a high-probability bound over $\X_n$. 
}

\label{tab:bq-rates}
\small
\setlength{\tabcolsep}{5pt}
\renewcommand{\arraystretch}{1.45}
\begin{adjustbox}{max width=\textwidth}
\begin{tabular}{lllll}
\toprule
Kernel
& Design 
& Noise Condition
& Sufficient $m$
& Error rate \\
\midrule

Mat\'ern
& Fixed
& Noiseless
& $n^{\frac{1}{s} \min\{s,\tau\}}$
& $n^{-\frac{2\min\{s, \tau\}}{d}}$ \\

Mat\'ern
& Random
& Noiseless
& $\left(\frac{n}{\log n} \right)^{\frac{1}{s} \min\{s,\tau\}}$
& $\left(\frac{\log n}{n}\right)^{\frac{2\min\{s, \tau\}}{d}}$ \\

Mat\'ern
& Fixed/Random
& Noisy
& $n^{d/(2s)}$
& $n^{-1}$ \\

Gaussian
& Fixed
& Noiseless
& $\frac{n}{(\log n)^{d}}$
& $n^{-\frac{2s}{d}} (\log n)^{2s}$ \\

Gaussian
& Random
& Noiseless
& $\frac{n}{(\log n)^{d+1}}$
& $n^{-\frac{2s}{d}} (\log n)^{2s(1+\frac{1}{d})}$ \\

Gaussian
& Fixed/Random
& Noisy
& $n^{d/(2s)}$
& $n^{-1}$ \\
\bottomrule
\end{tabular}
\end{adjustbox}
\end{table*}

Our second contribution is to show that the same operator-adapted features retain sharp statistical performance for related kernel methods such as kernel ridge regression (KRR), but with a different feature budget from that required for quadrature. For noisy observations, both the Mat\'ern and Gaussian constructions attain the Sobolev minimax squared population prediction rate
$
n^{-2s/(2s+d)}
$
under fixed and random designs. Compared with classical full-kernel KRR \cite{wang:2022}, our method requires only $m\asymp n^{d/(2s+d)}$ features, thus reducing computational complexity. Our finite-feature method also avoids saturation of full-kernel KRR when $s > 2\tau$ \cite{wang:2022, LiZhangLin:2023saturation}. In addition, we fill the theoretical gaps when the target does not belong to the kernel RKHS, and when Gaussian kernel is used for the present Sobolev/operator-feature setting. Our finite-feature framework extends naturally to kernel mean embedding, maximum mean discrepancy, and other tasks involving kernel integrals or kernel matrix inversion.

Our third contribution concerns uncertainty quantification. We study the frequentist coverage of the resulting credible interval. For $t\in (0,1)$, let $I_t$ denote the nominal $(1-t)$ credible interval centered at the finite-feature posterior mean. For quasi-uniform fixed designs, we prove the asymptotic coverage,
\begin{align*}
\liminf_{n\rightarrow\infty}
\inf_{\|f\|_{H^s(\Omega)}\leq 1}
P_\epsilon\left(
I(f)\in I_t\mid\X_n
\right)
\geq 1-t.
\end{align*}
As a result, the interval may be asymptotically conservative, but its coverage does not fall below the nominal level uniformly over the Sobolev unit ball. For i.i.d.\ random designs, concentration of the empirical feature Gram matrix additionally yields variance calibration and hence uniformly asymptotic coverage.

Finally, we identify concrete operator choices that make the feature integrals tractable without changing the statistical rates. 
Our theory shows that the convergence rate of finite-feature quadrature or regression error relies on the operator only through the asymptotic behavior of its eigenpairs.
Moreover, when $\Omega$ is contained in a box of radius $B$, our results require only a fixed padding of the operator domain $D_L$, for example with $L=B+1$. 
In fact, for the Laplacian/Legendre operator defined on $D_L$, allowing $L$ to diverge can lead to a loss in the convergence rate.
This differs from Hilbert space Gaussian process (HSGP) approximations that enlarge the auxiliary domain with sample size $n$ to improve global kernel or posterior approximation \cite{RiutortMayol:2023,Solin:2020}. In our examples,
Laplacian eigenfunctions yield Fourier-type features and are particularly suitable when $\PP$ has an analytically available characteristic function. Tensor-product Legendre eigenfunctions yield polynomial features and are suitable when the moments of $\PP$ are analytically available, which substantially extends previous works that have only focused on Fourier features. Our analysis can be extended to other differential operators whose eigenpairs are sufficiently tractable and satisfy the required ellipticity, regularity, spectral, and eigenfunction conditions.

\subsection{Related work and comparison}
\label{sec:related_work}

A substantial part of the kernel quadrature literature approaches integration through approximation of the kernel mean embedding. Let $\mathcal H_k$ denote the RKHS of kernel $k$.
For $f\in\mathcal H_k$, the integration error of a weighted rule can be bounded by
\begin{align*}
    \sup_{\|f\|_{\mathcal{H}_k} \leq 1} \left|\int_\Omega f(\bfx) \ud\PP(\bfx) - \sum_{i=1}^n w_i f(\bfx_i) \right|
    =
    \left\|\int_\Omega k(\cdot, \bfx) \ud\PP(\bfx)  - \sum_{i=1}^n w_i k(\cdot, \bfx_i) \right\|_{\mathcal{H}_k}.
\end{align*}
This relation connects quadrature error to kernel mean approximation and underlies random-feature and Nystr\"om approaches \cite{Bach:2017, Chatalic:2022}. Given a candidate pool of $N$ points, leverage-score sampling and related Nystr\"om procedures select $n$ evaluation points and construct a quadrature rule by projection or regularization \cite{Chatalic:2025,Chatalic:2022}. Other approaches use kernel recombination \cite{Adachi:2022BASQ, Hayakawa:2022KQ}, sampling-based Nystr\"om approximation \cite{Hayakawa:2023NystromKQ}, or kernel thinning \cite{Dwivedi:2024} to reduce function evaluations. These methods provide strong guarantees for RKHS worst-case error or kernel mean approximation, but differ from ours. They typically reduce a candidate pool or approximate the kernel mean, while we jointly replace the kernel matrix and kernel mean by an $m$-dimensional feature representation. Depending on the construction, their final rules may still require kernel algebra involving the selected $n$ points.

A second, closely related literature studies low-rank approximations for kernel ridge regression. Nystr\"om projection and computational regularization approximate the kernel matrix using sampled columns or inducing subspaces while retaining sharp regression rates \cite{Li:2023, Lin:2020, Rudi:2015, Rudi:2017FALKON}. Randomized kernel sketches provide statistically accurate low-dimensional systems \cite{ElAlaoui:2015, Yang:2017}, while random and structured Fourier features use explicit finite-dimensional feature maps \cite{Avron:2017, Li:2021, Rudi:2017, Doumeche:2026}. These methods establish important computational-statistical tradeoffs for regression, but their assumptions and losses do not directly resolve our quadrature problem. Much Nystr\"om theory formulates regularity through operator source conditions rather than Sobolev membership, while other results assume $f\in\mathcal H_k$ or primarily control empirical prediction loss. Structured Fourier analyses can provide population guarantees under more interpretable smoothness assumptions, but sharp results relevant here are generally restricted to specialized matched-smoothness or kernel regimes. Most importantly, low-rank regression guarantees do not control the separate approximation of $\int_\Omega k(\cdot,\bfx) \ud\PP(\bfx)$ and therefore do not imply an end-to-end quadrature bound.

Generalized spectral Bayesian quadrature and fast Fourier Bayesian quadrature construct rules from explicit spectral features, reducing matrix calculations and replacing the kernel mean by feature integrals \cite{Warren:2022, Warren:2024FFBQ}. Trigonometric quadrature Fourier features similarly use deterministic Fourier-domain quadrature for accurate kernel approximations in scalable Gaussian process regression \cite{Li:2024TQFF}. These methods demonstrate the computational value of spectral features and, in some cases, provide kernel-level approximation bounds. Such bounds, however, do not automatically imply estimator-level guarantees, since kernel error can be amplified by the regularized inverse, the feature approximation affects both the Gram matrix and kernel mean, and the integrand may lie outside the RKHS. 
Our construction uses ordered deterministic features instead of randomly sampled frequencies or Fourier-domain quadrature nodes, and analyzes the resulting quadrature estimator end to end.

Our feature representation is also structurally related to two deterministic Fourier feature Gaussian process approximations, the Hilbert space Gaussian process and the equispaced Fourier Gaussian process (EFGP). HSGP approximates a stationary kernel by a truncated expansion in Laplacian eigenfunctions on an auxiliary domain, weighted by the kernel spectral density \cite{Solin:2020, RiutortMayol:2023}, and has proved effective for reduced-rank GP regression and probabilistic programming. EFGP instead uses regular frequency grids and fast transforms for uniform kernel approximation and scalable GP regression \cite{Barnett:2024, Greengard:2025}. These works primarily analyze the kernel, covariance matrix, or GP posterior, but not the downstream quadrature error when both the regularized inverse and kernel integral are replaced. The auxiliary domain also plays a different role. For example, HSGP enlarges it to reduce boundary-induced kernel approximation error, while our estimator-level analysis permits a fixed separation between $\Omega$ and $\partial D_L$ and shows that unnecessarily increasing $L$ can reduce the effective spectral cutoff.

Another line of work approximates the Gaussian process posterior using $m < n$ inducing variables, reducing the cost of exact GP inference to approximately $nm^2 + m^3$. Existing theory evaluates these approximations under several criteria. \cite{Burt:2020SVGP} bound the reverse Kullback--Leibler (KL) divergence between the sparse variational and exact GP posteriors, closely related to the gap between the log marginal likelihood and evidence lower bound. \cite{Nieman:2022SVGP} establish frequentist contraction rates around the true regression function, while subsequent work studies global or pointwise variational credible sets \cite{Nieman:2023SVGPUQ, Travis:2023SVGPUQ}. More recent work derives a strictly tighter evidence lower bound at the same computational order \cite{Bui:2025TighterSVGP}. These results provide guarantees for posterior approximation and uncertainty quantification, but do not address the kernel integration problem. Thus, posterior-level KL divergence, contraction, or coverage guarantees do not directly yield end-to-end kernel quadrature bounds when both the kernel matrix and kernel mean are approximated.

Table~\ref{tab:related_works} summarizes representative results from these different lines of work. Because the methods address different statistical tasks, the table is intended to compare their assumptions, computational objects, and available guarantees rather than to give a strict numerical ranking. Relative to the cited literature, the distinctive scope of our analysis is the simultaneous treatment of an explicitly integrable finite-feature kernel, the regularized quadrature estimator, Sobolev targets that need not belong to the kernel RKHS, and coupled tuning of the number of features $m$, the regularizer/nugget $\eta$, and the lengthscale $\ell$, under general smoothness misspecification.

\begin{table*}[t]
\centering
\caption{Comparison of related kernel methods with random design. The guarantee is non-squared (worst-case) prediction error. Here we consider a Mat\'ern kernel with smoothness $\nu$ and $\tau = \nu+d/2$, $N$ = number of sample pool points; $n$ = number of function evaluations; $m$ = number of Nystr\"om landmarks or finite features; $n_{iter}$ = number of iterations in the conjugate gradient algorithm; $\tilde{O}(\cdot)$ is equal to $O(\cdot)$ up to logarithms; $f_H = \arg\min_{f^* \in \overline{\mathcal{H}_k}} \|f - f^* \|_{L^2}$, where $\overline{\mathcal{H}_k}$ is the closure of $\mathcal{H}_k$; for any $L^2$-integrable function $f$ with respect to $\PP$,
$\mathcal L_k(f)(\cdot) = \int_{\Omega} k(\bfx,\cdot)f(\bfx)\ud \PP(\bfx)$ is the kernel convolution operator.
}
\label{tab:related_works}
\small
\setlength{\tabcolsep}{5pt}
\renewcommand{\arraystretch}{1.45}
\begin{adjustbox}{max width=\textwidth}
\begin{tabular}{llll}
\toprule
Methods 
& Complexity
& Selected assumptions
& Guarantee \\
\midrule

\multicolumn{3}{l}{\textit{A. Quadrature and kernel-mean methods}}
& Noiseless quadrature error \\
\midrule

Nystr\"om leverage score sampling \cite{Chatalic:2025} 
&
$Nn^2 + n^3$
&
$f \in \mathcal{H}_k$; $N \gtrsim \tilde{O}(n^{2\tau/d})$
&
$\tilde{O}(n^{-\tau/d})$
\\

Kernel thinning \cite{Dwivedi:2024}
& $N^2$
& $N = O(n^2)$
& $\tilde{O}(n^{-1})$ \\

Nystr\"om kernel approximation \cite{Hayakawa:2023NystromKQ}
& $nmN + nm^2 + n^3 \log(N/n)$
& $f \in \mathcal{H}_k$; $N\gtrsim O(n^{2\nu/d})$, $m \gtrsim\tilde{O}(n^{2+\nu/d})$
& $O(n^{-\nu/d})$ \\

Fourier feature BQ \cite{Warren:2024FFBQ}
& $nm^2 + m^3$
& Stationary kernel
& No end-to-end theory \\

Our method
& $nm^2 + m^3$
& $s > d/2$
& $\tilde{O}(n^{-\min\{s, \tau\}/d})$ \\
\midrule

\multicolumn{3}{l}{\textit{B. KRR and GP approximations}}
& Regression error \\
\midrule

Nystr\"om projection KRR \cite{Lin:2020} 
& $nm^2 + m^3$
& $f \in L^\infty$; $f_H = \mathcal{L}_k^{\zeta} g, \zeta\ge 0, \int |g|^2 \ud\PP < \infty$
& $\tilde{O}\left(n^{-\frac{\min\{\zeta, 1\}}{2\min\{\zeta, 1\} + d/(2\tau)}} \right)$ \\

Kernel sketching \cite{Yang:2017}
& $n^2\log m + m^3$
& $f \in \mathcal{H}_k$; $m \gtrsim O(n^{\frac{d}{2\tau +d}})$
& Empirical error: $O(n^{-\frac{\tau}{2\tau + d}})$ \\

Random feature KRR \cite{Li:2021}
& $nm^2 + m^3$
& $m \gtrsim \tilde{O}(n^{d/(4\tau)})$
& $O(n^{-1/4})$ \\

Fourier feature KRR \cite{Doumeche:2026} 
& $n\log n + m^2 \log m$
& $s = \tau \in \NN$; $m \gtrsim O(n^{\frac{d}{2\tau+d}})$
& $O(n^{-\frac{\tau}{2\tau + d}})$ \\

Our method
& $nm^2 + m^3$
& $s > d/2$
& $O(n^{-\frac{s}{2s+d}})$ \\
\bottomrule

Equispaced Fourier GP (EFGP) \cite{Greengard:2025}
& $n + n_{iter} m \log m$
& Stationary kernel
& Kernel-level guarantee \\

Hilbert Space GP (HSGP) \cite{Solin:2020} 
& $nm^2 + m^3$
& Stationary kernel
& Kernel-level guarantee \\
\bottomrule

\end{tabular}
\end{adjustbox}
\end{table*}

We define some universal notation. For nonnegative quantities $a$ and $b$, we write $a\lesssim b$ if $a\leq Cb$ for a constant $C>0$, $a\gtrsim b$ if $b\lesssim a$, and $a\asymp b$ if both $a\lesssim b$ and $b\lesssim a$. 
The value of $C$ may change from line to line. Unless stated otherwise, implicit constants may depend on the fixed structural quantities in the relevant result, including the dimension $d$ and smoothness parameters $s, \nu$, fixed kernel parameters $\sigma^2, \ell$, the domains $\Omega$, density and mesh-ratio bounds, and the coefficient, regularity, and ellipticity constants of the differential operator $\mathcal T_L$, but not on $n,m,\eta,\lambda,\delta,\sigma_\epsilon$, or on $\ell$ and $L$ when these are allowed to vary and their dependence is displayed explicitly. 
We write $U\Subset V$ if $\overline U$ is a compact subset of $V$, and write $a_+:=\max\{a,0\}$. 
For a generic vector $\bfa$, $\|\bfa\|_2$ denotes its Euclidean norm. For a generic function $g$, let $\|g\|_{L^q(\Omega)}=\left\{\int_\Omega |g(\bfx)|^q \ud \bfx \right\}^{1/q}$ and $\|g\|_n=\left\{n^{-1}\sum_{i=1}^n g^2(\bfx_i)\right\}^{1/2}$. Domains or measures are indicated explicitly when needed.
The symbols $\EE_\epsilon$ and $P_\epsilon$ refer to expectation and probability over the observation noise conditionally on $\X_n$.
We use $\NN=\{0,1,2,\ldots\}$ for the set of all nonnegative integers and $\NN_+=\{1,2,\ldots\}$ for the set of all positive integers.

The remainder of the paper is organized as follows.
Section~\ref{sec:operator_approximation} introduces the operator-adapted finite-feature construction and the resulting kernel quadrature estimator.
Section~\ref{sec:main_theory} establishes non-asymptotic error bounds for Mat\'ern and Gaussian kernels.
Section~\ref{sec:posterior_variance} studies the frequentist coverage of the finite-feature confidence intervals.
Section~\ref{sec:concrete_bases} specializes the general theory to several concrete examples of operators and derives explicit tuning rules.
Section~\ref{sec:feature_lower_bound} analyzes the lower bound of feature budget for different tasks.
Section~\ref{sec:related_estimators}
discusses consequences for related kernel estimators.
Numerical results are presented in Section~\ref{sec:numerical_experiments},
and Section~\ref{sec:discussion} concludes the paper.

\section{Operator-adapted finite-feature kernel quadrature estimation}
\label{sec:operator_approximation}
We first introduce the operator-adapted kernel estimator for quadratures based on finite features. The main idea is to replace the generally intractable eigenfunctions in the Mercer expansion of a kernel by the tractable eigenfunctions of a differential operator defined on a slightly larger domain.

\subsection{Operator-based kernel approximation}
We illustrate the main idea using the isotropic Mat\'ern kernel \cite{Rasmussen:2006book, Stein:1999book} as a working example, which is defined as
\begin{align} \label{eq:kernel_matern}
    k(\bfx, \bfx^\prime) = \sigma^2 \frac{2^{1-\nu}}{\Gamma(\nu)} \left(\frac{\sqrt{2\nu} \|\bfx-\bfx^\prime\|_2}{\ell}\right)^\nu K_\nu\left(\frac{\sqrt{2\nu} \|\bfx-\bfx^\prime \|_2}{\ell}\right),~~ \bfx,\bfx^\prime \in \RR^d,
\end{align}
where $K_\nu$ is the modified Bessel function of the second kind, $\ell$ is the lengthscale, $\sigma^2$ is the scale, and $\nu > 0$ is the smoothness. 

Under the usual conditions for Mercer expansion, 
the kernel \eqref{eq:kernel_matern} admits the Mercer representation $k(\bfx,\bfx') = \sum_{j=1}^\infty \mu_j\psi_j(\bfx)\psi_j(\bfx')$. Truncating this expansion would give a natural finite-feature approximation. However, for the Mat\'ern kernels, the eigenvalues $\mu_j$, eigenfunctions $\psi_j$ and their integrals with respect to $\PP$ generally do not have closed forms except in a few special cases \cite{vantrees:2001book}, limiting their applicability to computationally efficient kernel quadrature.
Instead of using the Mercer expansion, we seek an alternative approach based on the kernel's spectral properties.
Let $\mathcal C_k:L^2(\RR^d)\to L^2(\RR^d)$ denote the whole-space convolution operator associated with the Mat\'ern kernel, $\mathcal C_k f(\bfx) = \int_{\RR^d} k(\bfx-\bfy) f(\bfy) \ud \bfy$, for any $f\in L^2(\RR^d)$. Since the Mat\'ern kernel belongs to $L^1(\RR^d)$, $\mathcal C_k$ is a bounded operator on $L^2(\RR^d)$. Under the convention of Fourier transform
$
    (\mathcal F f)(\bm\omega) = \int_{\RR^d} \exp\left(-i \bm\omega^\top \bfx\right) f(\bfx) \ud \bfx
$ where $i=\sqrt{-1}$,
$\mathcal C_k$ acts by multiplication with the Matérn spectral density $S$:
\begin{align*}
    (\mathcal F (\mathcal C_k f))(\bm\omega) 
    &= (\mathcal F k)(\bm\omega) ~(\mathcal F f)(\bm\omega) 
    \\
    &
    \propto \left(\frac{2\nu}{\ell^2}+\|\bm\omega\|_2^2\right)^{-\tau} ~(\mathcal F f)(\bm\omega) 
    = \mathcal F \left(\left(\frac{2\nu}{\ell^2} - \Delta \right)^{-\tau} f \right) (\bm\omega),
\end{align*}
and it follows that
\begin{align} \label{eq:convolution_op_matern}
    \mathcal C_k \propto \left(\frac{2\nu}{\ell^2} - \Delta \right)^{-\tau}.
\end{align}
Thus, Fourier features arise naturally because they diagonalize the Laplacian operator. In \eqref{eq:convolution_op_matern}, $\left(\frac{2\nu}{\ell^2} - \Delta \right)^{-\tau}$ is also known as the Bessel operator \cite{Secchi:2017}. 

The representation \eqref{eq:convolution_op_matern} motivates replacing the Fourier basis by eigenfunctions of a tractable differential operator on a bounded auxiliary domain. We introduce a second-order differential operator $\mathcal{T}_L$ on an auxiliary bounded domain $D_L \Supset \Omega$, and replace $-\Delta$ in \eqref{eq:convolution_op_matern} by $\mathcal T_L$. We then define the operator-induced covariance operator
$
    \left(\frac{2\nu}{\ell^2} + \mathcal T_L \right)^{-\tau}.
$
If $\mathcal{T}_L$ is non-negative and self-adjoint with compact resolvent on $L^2(D_L)$, then it admits an eigensystem
\begin{align} \label{eq:T_L.eigen}
    \mathcal{T}_L \phi_j &= \lambda_j \phi_j, ~j \in \NN_+,
\end{align}
where the eigenvalues are $0 \leq \lambda_1 \leq \lambda_2 \leq \ldots$ and the eigenfunctions
$\{\phi_j\}_{j \in \NN_+}$ form an orthonormal basis in $L^2(D_L)$. Spectral functional calculus then gives
\begin{align*}
    \left(\frac{2\nu}{\ell^2} + \mathcal T_L \right)^{-\tau} \phi_j = \left(\frac{2\nu}{\ell^2} + \lambda_j \right)^{-\tau} \phi_j.
\end{align*}
Therefore, the following finite-feature sum is a working kernel motivated by the Mat\'ern spectral profile:
\begin{align} \label{eq:kernel_approx_matern}
    \widehat k_m(\bfx, \bfx') := \sum_{j = 1}^m \left(\frac{2\nu}{\ell^2} + \lambda_j \right)^{-\tau} \phi_j(\bfx) \phi_j(\bfx'), ~~ \bfx,\bfx^\prime \in \Omega.
\end{align}
More generally, the same construction applies to other covariance classes by replacing the weights $\left(\frac{2\nu}{\ell^2} + \lambda_j \right)^{-\tau}$ with an appropriate nonnegative spectral weight. In particular, we will later consider Gaussian-type weights with exponential spectral decay.

It is useful to distinguish the construction in \eqref{eq:kernel_approx_matern} from a truncated Mercer expansion. The eigenfunctions $\phi_j$ are determined by the differential operator $\mathcal T_L$, its boundary conditions, and the auxiliary domain $D_L$,  while they are constructed independently of the integration measure $\PP$. After restriction from $D_L$ to $\Omega$, the functions $\{\phi_j\}$ need not be orthogonal in $L^2(\Omega)$ and need not coincide with the Mercer eigenfunctions $\{\psi_j\}$. This distinction is important for kernel quadrature because the feature integrals $\int_\Omega
\phi_j(\bfx)\ud\PP(\bfx)$ can remain analytically tractable even when the corresponding Mercer eigenfunctions are analytically intractable.

For a general differential operator $\mathcal T_L$, we also emphasize that we do not require the operator-adapted kernel $\widehat k_m$ to be a uniformly accurate approximation to the original Mat\'ern kernel. Rather, the Mat\'ern spectral representation motivates the polynomially decaying weights $\left(\frac{2\nu}{\ell^2} + \lambda_j \right)^{-\tau}$, and the statistical properties of the resulting finite-feature estimator are analyzed directly. This perspective differs from approaches whose primary objective is to approximate the kernel itself. In our setting, the relevant question is whether the finite-feature representation preserves the convergence rate of the resulting quadrature estimator. 

Unlike operator-stochastic partial differential equation (SPDE) constructions \cite{Lindgren:2011SPDE, Cox:2020, Borovitskiy:2020RiemannianMatern} in which the Gaussian field and the elliptic operator are defined on the same spatial domain $\Omega$, we define the operator on an auxiliary domain $D_L$ and restrict its eigenfunctions to an interior target domain $\Omega \Subset D_L$. This separation permits the use of interior elliptic regularity and avoids imposing artificial boundary conditions directly on the integration domain. Consequently, the auxiliary domain structure enables an end-to-end error analysis for the resulting quadrature estimator when both the regularized Gram matrix and the kernel mean are replaced by the finite-feature representations. More importantly, certain choices of the operator $\mathcal T_L$ together with boundary conditions lead to closed-form eigenpairs $\{(\lambda_j,\phi_j),~j\in \NN_+\}$ and integrals $\int_\Omega
\phi_j(\bfx)\ud\PP(\bfx)$, which enables efficient computation for kernel quadrature estimation. Concrete constructions based on Laplacian and tensor-product Legendre eigenfunctions are presented in Section~\ref{sec:concrete_bases}.

\subsection{Finite-feature predictor for quadrature}

We consider the observation model
\begin{align*}
    y(\bfx_i)
    &=
    f(\bfx_i)+\epsilon_i, 
    \quad 
    \epsilon_i \sim \mathcal{N}(0, \sigma_\epsilon^2), 
    \qquad
    i=1,\ldots,n, 
\end{align*}
where $\bfx_i\in\Omega=(-B,B)^d
$, $f:\Omega\rightarrow\mathbb R$ is the unknown integrand, $\epsilon_i$ denotes observation noise, and $\sigma_\epsilon^2 \geq 0$ is a constant. The noiseless setting is recovered by taking $\epsilon_i=0$. Let
\begin{align*}
    \X_n
    :=
    \{\bfx_1,\ldots,\bfx_n\},
    \qquad
    \Y_n
    :=
    \left[
    y(\bfx_1),\ldots,y(\bfx_n)
    \right]^\top.
\end{align*}
Our objective is to estimate
$
    I(f)
    :=
    \int_\Omega f(\bfx) \ud\PP(\bfx),
$
where $\PP$ is a probability measure on $\Omega$.

For the Bayesian quadrature construction, we assign a Gaussian process prior with kernel $k(\cdot,\cdot)$,
$
    f
    \sim
    \mathcal{GP}(0,k).
$
The posterior distribution of $I(f)$ is Gaussian, with mean and variance
\begin{align*}
    \EE\left[
    I(f)\mid\Y_n
    \right]
    &=
    \bm\mu_n^\top
    \left(
    \K_n+\sigma_\epsilon^2\I_n
    \right)^{-1}
    \Y_n,
    \\
    \Var\left[
    I(f)\mid\Y_n
    \right]
    &=
    \iint_{\Omega \times \Omega} k(\bfx, \bfx^\prime) \ud\PP(\bfx) \ud\PP(\bfx^\prime)
    -
    \bm\mu_n^\top
    \left(
    \K_n+\sigma_\epsilon^2\I_n
    \right)^{-1}
    \bm\mu_n,
\end{align*}
where $\bm\mu_n = \left[ \int_\Omega k(\bfx_1, \bfx) ~\ud\PP(\bfx),\ldots, \int_\Omega k(\bfx_n,\bfx) ~\ud\PP(\bfx) \right]^\top$, $\K_n=\{k(\bfx_i,\bfx_j)\}_{1\le i,j\le n}$ is the $n\times n$ kernel matrix, and $\I_n$ is the $n\times n$ identity matrix.

Motivated by the form of exact Bayesian quadrature, we define a separate finite-dimensional ``working'' Gaussian process, using the finite-feature kernel in \eqref{eq:kernel_approx_matern}:
\begin{align*}
    \bar f_m \sim \mathcal{GP}(0, \sigma_\epsilon^2 \eta^{-1} \widehat k_m), 
    \qquad \widehat k_m(\bfx, \bfx^\prime) = \sum_{j=1}^m w(\lambda_j) \phi_j(\bfx) \phi_j(\bfx^\prime), 
    \qquad
    \epsilon_i \sim \mathcal{N}(0, \sigma_\epsilon^2),
\end{align*}
where $\eta>0$ is a tunable nugget, $\{(\lambda_j,\phi_j):j\in \NN_+\}$ are the eigenpairs of $\mathcal{T}_L$ in \eqref{eq:T_L.eigen}, and the weights $w(\lambda_j)$ are chosen later according to the spectral decay of kernel $k$. Let
\begin{align*}
    \bm\Phi=(\phi_j(\bfx_i))_{1\le i\le n,\,1\le j\le m},
    \quad
    \W=\operatorname{diag}(w(\lambda_1),\ldots,w(\lambda_m)),
    \quad
    \widehat \K_n= \bm\Phi 
    \W \bm\Phi^\top.
\end{align*}
The working posterior distribution is Gaussian with mean and variance:
\begin{align} \label{eq:working_model_dist}
    \EE\left[
    I(\bar f_m)\mid\Y_n
    \right]
    &=
    \widehat{\bm\mu}_n^\top
    \left(
    \widehat \K_n+\eta\I_n
    \right)^{-1}
    \Y_n,
    \nonumber \\
    \Var\left[
    I(\bar f_m)\mid\Y_n
    \right]
    &=
    \frac{\sigma_\epsilon^2}{\eta} \left[\iint_{\Omega \times \Omega} \widehat k_m(\bfx, \bfx^\prime) \ud\PP(\bfx) \ud\PP(\bfx^\prime)
    -
    \widehat{\bm\mu}_n^\top
    \left(
    \widehat\K_n+\eta\I_n
    \right)^{-1}
    \widehat{\bm\mu}_n \right],
\end{align}
where $\widehat{\bm\mu}_n = \left[ \int_\Omega \widehat k_m(\bfx_1, \bfx) ~\ud\PP(\bfx),\ldots, \int_\Omega \widehat k_m(\bfx_n,\bfx) ~\ud\PP(\bfx) \right]^\top$.
Let $\bm\phi(\bfx) = \left[\phi_1(\bfx), \ldots, \phi_m(\bfx) \right]^\top$, then $\widehat{\bm\mu}_n = \bm\Phi \W \left[\int_\Omega \bm\phi(\bfx) \ud \PP(\bfx) \right]$.
The resulting finite-feature quadrature estimator is
\begin{align} \label{eq:finite_feature_estimator}
    \widehat I_m(f)
    &=
    \EE\left[
    I(\bar f_m)\mid\Y_n
    \right]
    =
    \widehat{\bm\mu}_n^\top
    \left(
    \widehat \K_n+\eta\I_n
    \right)^{-1}
    \Y_n .
\end{align}
In the noiseless case, $\sigma_\epsilon^2=0$, so $\Y_n=\F_n$, while $\eta>0$ is retained for regularization.
Since
\begin{align*}
    \left(
    \widehat \K_n+\eta\I_n
    \right)^{-1} \widehat{\bm\mu}_n
    =
    \bm\Phi \left( \bm\Phi^\top \bm\Phi + \eta \W^{-1} \right)^{-1} \left[\int_\Omega \bm\phi(\bfx) \ud \PP(\bfx) \right],
\end{align*}
computing $\widehat I_m(f)$ reduces to solving an $m$-dimensional linear system with coefficient matrix $\bm\Phi^\top\bm\Phi+\eta\W^{-1}$. Thus, the original $O(n^3)$ computational complexity for the BQ posterior mean and variance is reduced to $O(nm^2 + m^3)$. This is beneficial especially when $m \lesssim n$. For Fourier features, this complexity may be further reduced by numerical algorithms \cite{Warren:2024FFBQ, Doumeche:2026}.

We consider two principal choices of spectral weights. For the Mat\'ern-type weights,
\begin{align} \label{eq:weights_matern}
    w(\lambda_j) = \left(\frac{2\nu}{\ell^2}+\lambda_j\right)^{-\tau}, \quad \tau = \nu+d/2.
\end{align}
For the Gaussian-type weights,
\begin{align} \label{eq:weights_gaussian}
    w(\lambda_j) = \ell^{d} \exp \left(-\frac{1}{2} \ell^2 \lambda_j\right).
\end{align}

Equations \eqref{eq:weights_matern} and \eqref{eq:weights_gaussian} specify Mat\'ern-type and Gaussian-type spectral profiles for an operator-induced working kernel. For a fixed $L$ and a general realization of $\mathcal{T}_L$, we do not attempt to show that 
$\widehat k_m$ converges to the original stationary kernel as $m \rightarrow \infty$. The infinite-series limit is instead the kernel associated with $w(\lambda_j)$ and $\phi_j$, which depends on the differential operator and its boundary realization. Our objective is to control the downstream quadrature error over Sobolev classes, rather than the approximation error to a prescribed stationary kernel. We also omit the scaling parameter $\sigma^2$ and the multiplicative constant in the Mat\'ern kernel \eqref{eq:kernel_matern}, since the scaling constant can be absorbed in the nugget $\eta$.

\section{Convergence analysis for quadrature}
\label{sec:main_theory}
We now establish convergence guarantees for the operator-adapted kernel quadrature estimator.
We measure the performance of the finite-feature estimator by the conditional worst-case mean squared error
\begin{align*} 
    \mathcal Q_n
    :=
    \sup_{\|f\|_{H^s(\Omega)}\leq1}
    \EE_\epsilon
    \left[
        \left(
            I(f)-\widehat I_m(f)
        \right)^2
        ~\middle|~
        \X_n
    \right].
\end{align*}
This criterion accounts for both the deterministic quadrature error and the contribution of observation noise, while treating the design $\X_n$ as fixed. For random designs, the resulting bound will additionally hold with high probability over the draw of $\X_n$. 

\subsection{Assumptions}
We first introduce the assumptions needed for our analysis.

\begin{assumption} \label{assump:noise_vector}
    Let $\bm\epsilon_n=\left[\epsilon_1, \ldots, \epsilon_n \right]^\top$ denote the noise vector. Suppose $\bm\epsilon_n \mid \X_n \sim \mathcal{N}(\bm 0, \sigma_\epsilon^2 \I_n)$. 
\end{assumption}

\begin{assumption} \label{assump:finite_density}
    The probability measure $\PP$ on $\Omega=(-B,B)^d$ admits a density $p(\cdot)$ with respect to the Lebesgue measure such that $p \in L^\infty(\Omega)$.
\end{assumption}

\begin{assumption} \label{assump:bounded_density}
    There exist constants $p_{\min}, p_{\max}$ such that $0 < p_{\min} \leq p(\bfx) \leq p_{\max} < \infty$ for all $\bfx\in \Omega$.
\end{assumption}

\begin{assumption} \label{assump:f_smoothness}
    The underlying function $f$ satisfies $f \in H^s(\Omega)$ with  $s > d/2$.
\end{assumption}

Assumption~\ref{assump:noise_vector} controls the variance introduced by the noisy observations. Gaussianity is not essential for the mean-squared error bounds in Section~\ref{sec:main_theory}, which only use the first two conditional moments of the noise, but it is used in Section~\ref{sec:posterior_variance} for the frequentist coverage analysis. Assumption~\ref{assump:finite_density} allows integration with respect to $\PP$ to be controlled through Lebesgue $L^2(\Omega)$ norms. For random designs, Assumption~\ref{assump:bounded_density} additionally provides a uniform comparison between the population design measure and the Lebesgue measure. The lower bound is used to control the population feature Gram matrix.
The condition $s>d/2$ in Assumption~\ref{assump:f_smoothness} implies that every element of $H^s(\Omega)$ has a continuous and bounded representative. In particular, the point evaluations $f(\bfx_i)$ appearing in the estimator are well defined, and Sobolev embedding inequalities can be applied.

Next, we introduce the assumptions on the differential operator.
Let $D_L := (-L, L)^d$ with $L>B$ such that $D_L\Supset \Omega$, and consider the second-order differential expression 
\begin{align} \label{eq:second_order_operator}
    T_L 
    = 
    -\sum_{i,j=1}^d \partial_i (A_{ij} \partial_j) + V.
\end{align}
Such divergence-form operators are standard in elliptic operator theory; see Chapters 4--5 of \cite{Ouhabaz:2005} and Equation 1.2 of \cite{Cox:2020}.

\begin{assumption} \label{assump:domain_padding}
    There exists a constant $c_0 > 0$ independent of $L$, such that
    $
        L \geq B + c_0.
    $
\end{assumption}

\begin{assumption} \label{assump:operator_coeff_smoothness}
    For the divergence form~\eqref{eq:second_order_operator}, there exist constants $M > 0$ and $Q \in \NN_+, Q > \max\{s, \tau\}/2$ such that 
    \begin{align*} 
        \sup_{L\geq B+c_0} \left(\|A_{ij}\|_{W_\infty^{2Q}(D_L)}+\|V\|_{W_\infty^{2Q}(D_L)} \right)
        \leq 
        M.
    \end{align*}
    Moreover, $\A = [A_{ij}]_{1\leq i,j\leq d}$ is symmetric, i.e., $A_{ij} = A_{ji}$.
\end{assumption}

\begin{assumption} \label{assump:operator_coeff_smoothness_gaussian}
    The conditions in Assumption~\ref{assump:operator_coeff_smoothness} hold, except that the regularity requirement on $Q$ is replaced by $ Q>s/2$.
\end{assumption}

\begin{assumption} \label{assump:local_elliptic}
    The differential expression $T_L$ in \eqref{eq:second_order_operator} is locally uniformly elliptic in $D_L$ (see Section 2, \cite{Metafune:2012, Versano:2023}). That is, for every fixed open set $U\Subset D_L$, there exists a constant $a_0 > 0$, dependent on $U, L$, such that
    \begin{align*}
        \bm\xi^\top \A(\bfx) \bm\xi
        \geq
        a_0 \|\bm\xi\|_2^2
    \end{align*}
    for every $\bm\xi\in\RR^d$ and almost every $\bfx\in U$. In addition, if $L \geq B + c_0$, and $U \Subset D_{B+c_0} = (-B-c_0, B+c_0)^d$ is independent of $L$, then $a_0$ depends only on $U, c_0, d$.
\end{assumption}

\begin{assumption} \label{assump:eigenpairs}
    The differential expression $T_L$ admits a nonnegative self-adjoint realization $\mathcal T_L$ on $L^2(D_L)$ with compact resolvent. We denote its domain by $\mathcal D(\mathcal T_L)$ and assume that
    $C_c^\infty(D_L)\subset\mathcal D(\mathcal T_L)$.
    By the spectral theorem for self-adjoint operators with compact resolvent (Theorem 6.11 in \cite{Brezis:2011}), $\mathcal{T}_L$ admits an eigensystem,
    \begin{align*} 
        \mathcal{T}_L \phi_j = \lambda_j \phi_j,
    \end{align*}
    where $\{\phi_j\}_{j \in \NN_+}$ is a complete orthonormal basis on $L^2(D_L)$, and $\{\lambda_j\}_{j\in \NN_+}$ are sorted in ascending order $0 \leq \lambda_1 \leq \lambda_2 \leq \ldots$.
\end{assumption}

\begin{assumption} \label{assump:eigenfunction_derivative}
    Each real partial derivative of $\phi_j$ exists
    on $\Omega$ and extends continuously to $[-B,B]^d$.
    There is a constant $C\ge1$, independent of $m$ and $\bm\alpha$, such that for every $m\ge1$ and every
    $\bm\alpha=(\alpha_1,\ldots,\alpha_d)\in\NN^d$,
    \begin{align*}
        \sup_{\bfx \in [-B,B]^d}
        \sum_{j=1}^m |D^{\bm\alpha}\phi_j(\bfx)|^2
        \le 
        C^{2|\bm\alpha|+2} (1+\lambda_m)^{|\bm\alpha|+d/2}.
    \end{align*}
    Here we use the usual multi-index notation $|\bm\alpha|=\sum_{i=1}^d \alpha_i$ and
    $D^{\bm\alpha}=\partial_{x_1}^{\alpha_1}\cdots
    \partial_{x_d}^{\alpha_d}$.
\end{assumption}

Assumptions~\ref{assump:domain_padding}--\ref{assump:eigenfunction_derivative} specify the analytic and spectral properties of the differential operator's domain and eigenpairs. Assumption~\ref{assump:domain_padding} ensures that the target domain $\Omega$ remains uniformly separated from the artificial boundary of $D_L$.
This buffer region allows functions on $\Omega$ to be extended and localized within $D_L$, while keeping the constants in the interior estimates uniform in $L$. 
Assumption~\ref{assump:operator_coeff_smoothness} is used for Mat\'ern-type weights. The condition
$
    2Q>\max\{s,\tau\}
$
provides sufficient regularity to control both the target smoothness and the Mat\'ern feature-space norm.
Assumption~\ref{assump:operator_coeff_smoothness_gaussian} is for Gaussian weights. 
Assumption~\ref{assump:local_elliptic} provides local uniform ellipticity on interior subsets of $D_L$. Together with Assumption~\ref{assump:operator_coeff_smoothness} or \ref{assump:operator_coeff_smoothness_gaussian}, it allows the Sobolev regularity of finite eigenfunction expansions on $\Omega$ to be controlled through spectral powers of $\mathcal T_L$. The interior ellipticity is required to accommodate operators such as the Legendre operator whose leading coefficients may degenerate at $\partial D_L$.
Assumption~\ref{assump:eigenpairs} guarantees the existence of the eigenpairs used in the finite-feature expansion. 
Different operator realizations and boundary conditions may produce different eigenvalues and eigenfunctions, but the theory depends only on the stated spectral and interior regularity properties.
Assumption~\ref{assump:eigenfunction_derivative} controls the growth of eigenfunction derivatives.
When $\bm\alpha=\bm{0}$, it is the spectral-sum bound
$\sum_{j=1}^m |\phi_j(\bfx)|^2 \leq C^2 (1+\lambda_m)^{d/2}$.
The same $C$ must work for all derivative orders.
It may depend on $L$, $B$, $d$, and eigensystem, but is independent of $m$. If $L = B+c_0$ is a fixed constant, then $C$ is also a fixed constant.

Finally, we state the assumptions on the design points.
Let
\begin{align*}
    h_n
    :=
    \sup_{\bfx\in\Omega}
    \min_{1\leq i\leq n}
    \|\bfx-\bfx_i\|_2,  \qquad
    q_n
    :=
    \frac12
    \min_{i\neq j}
    \|\bfx_i-\bfx_j\|_2 ,
\end{align*}
where $h_n$ denotes the fill distance of $\{\bfx_1,\ldots,\bfx_n\}$, and $q_n$ denotes their separation radius.

\begin{assumption} \label{assump:quasi_uniform}
    The design $\X_n$ is fixed and its mesh ratio is uniformly bounded:
    $
        \sup_{n \geq 2} ~h_n/q_n < \infty.
    $
\end{assumption}

\begin{assumption} \label{assump:random_iid}
    The design $\X_n$ is random and $\{\bfx_1, \ldots, \bfx_n\} \overset{\mathrm{i.i.d.}}{\sim} \PP$. 
\end{assumption} 

Assumption~\ref{assump:quasi_uniform} permits the use of deterministic Sobolev sampling inequalities and converts bounds in $h_n$ into explicit rates in $n$. It is also known as the quasi-uniform condition and indicates $h_n \asymp q_n \asymp n^{-1/d}$. 
Under Assumption~\ref{assump:random_iid}, deterministic sampling arguments are replaced by concentration inequalities for the empirical feature Gram matrix and the empirical approximation error. 

\subsection{Convergence rate of worst case error with Mat\'ern-type weights}
\label{subsec:matern_theory}
Suppose $w(\lambda_j)$ is the Mat\'ern weight function defined in \eqref{eq:weights_matern}. We first establish non-asymptotic error bounds for general operator-adapted feature systems. These results do not require the target function to belong to the Mat\'ern RKHS. Instead, the target smoothness $s$ and the kernel regularity $\tau=\nu+d/2$ are allowed to differ, covering both the oversmoothing and undersmoothing regimes. A notable feature of the following bounds is that all principal approximation and statistical parameters remain explicit. In particular, the spectral cutoff enters through the $m$th smallest eigenvalue $\lambda_m$, the sampling geometry through $h_n$, and regularization through the nugget $\eta$. Consequently, the bounds can be used not only to establish convergence, but also to derive concrete tuning rules for $m$ and $\eta$.

We first consider a deterministic quasi-uniform design. The resulting bound separates the spectral approximation error and the contribution of observation noise.

\begin{theorem}\label{thm:bq_matern_fixed} 
Suppose that Assumptions~\ref{assump:noise_vector},
\ref{assump:finite_density},
\ref{assump:f_smoothness},
\ref{assump:domain_padding},
\ref{assump:operator_coeff_smoothness},
\ref{assump:local_elliptic},
\ref{assump:eigenpairs}, and
\ref{assump:quasi_uniform} hold, and $w(\lambda_j)$ is the Mat\'ern weight function in \eqref{eq:weights_matern}.
If $n$ is sufficiently large, then
\begin{align*} 
    \mathcal Q_n &\lesssim \Bigg\{ h_n^{d/2} \left( h_n^\nu+\eta^{1/2} \right) (1+\lambda_m)^{\frac{1}{2}(\tau-s)_+} 
    + \left( \eta^{-1/2}h_n^\nu+1 \right) \left[ h_n^s+(1+\lambda_m)^{-s/2} \right] \Bigg\}^2 \\ 
    &\qquad\quad + \sigma_\epsilon^2 \left( \eta^{-1/2}h_n^\tau+h_n^{d/2} \right)^2.
\end{align*} 
The hidden multiplicative constant is independent of $n$, $m$, $\eta$, and $\sigma_\epsilon^2$. 
\end{theorem} 

The upper bound in Theorem~\ref{thm:bq_matern_fixed}, as well as those in the remainder of this section and the next, is stated in an abstract form. We defer the explicit choices of the tuning quantities and the resulting convergence rates in $n$ to Section~\ref{sec:concrete_bases}. The upper bound for the quadrature error $\mathcal Q_n$ depends crucially on the feature dimension $m$ through $\lambda_m$, and on the regularizer/nugget $\eta$.  The factor
$
(1+\lambda_m)^{\frac{1}{2}(\tau-s)_+}
$
quantifies the smoothness mismatch between the Mat\'ern feature space and the target Sobolev space. This factor disappears when $s\geq\tau$, and produces an explicit oversmoothing penalty when $s<\tau$. 
The remaining terms describe the interaction between the design resolution and regularization. In particular, increasing $m$ reduces the truncation error, but may amplify the smoothness-mismatch term, while decreasing $\eta$ reduces regularization bias but enlarges the norm of the quadrature weights. Therefore, Theorem~\ref{thm:bq_matern_fixed} makes the relevant bias--variance trade-off explicit. 
For noisy observations, the final term isolates the effect of the true noise level $\sigma_\epsilon^2$. In the noiseless setting, the same result applies after setting $\sigma_\epsilon^2=0$.
Theorem~\ref{thm:bq_matern_fixed} also shows that the effects of the differential operator enter only through its eigenvalue growth and eigenfunction concentration, which permits the same analysis to accommodate different operator-adapted bases.

We next consider random design where the design points are independently sampled. In this setting, we provide a high-probability bound over the design $\X_n$.
For $\Lambda \geq 1$ and $m\in \NN_+$, let
\begin{align} \label{eq:N.H}
    N(\Lambda)
    :=
    \#\{j: 1+\lambda_j \leq \Lambda\},
    \qquad
    H(m)
    :=
    \sup_{\bfx\in\Omega}
    \sum_{j=1}^m
    \left|\phi_j(\bfx) \right|^2
\end{align}
Here $N(\Lambda)$ measures the spectral complexity of the differential operator, while $H(m)$ controls the pointwise concentration of the first $m$ eigenfunctions. The conditions on $n$ below are stated directly in terms of these two quantities and apply to general operator eigenpairs.

\begin{theorem}\label{thm:bq_matern_random}
Suppose that Assumptions~\ref{assump:noise_vector},
\ref{assump:bounded_density},
\ref{assump:f_smoothness},
\ref{assump:domain_padding},
\ref{assump:operator_coeff_smoothness},
\ref{assump:local_elliptic},
\ref{assump:eigenpairs}, and
\ref{assump:random_iid} hold, and $w(\lambda_j)$ is the Mat\'ern weight function in \eqref{eq:weights_matern}.
Let $R_m=1+\lambda_m$ and 
$
    L_p = \log \left( 2^{2p+4}N(2^{p+1}R_m) /\delta \right), ~p=0,1,\ldots.
$
For any $\delta\in(0,1)$, suppose 
\begin{align*} 
    n \gtrsim \left(\frac{\eta}{n}\right)^{-d/(2\tau)} \log\left(\frac{4m}{\delta}\right), 
    \qquad
    n \gtrsim \sum_{p=0}^{\infty} 2^{-sp} H\left(N(2^{p+1}R_m)\right) L_p,
\end{align*} 
where the implicit constant in the first lower bound on $n$
is sufficiently large.
Then, with probability at least $1-\delta$ with respect to $\X_n\overset{\mathrm{i.i.d.}}{\sim}\PP$, 
\begin{align*} 
    \mathcal{Q}_n
    &\lesssim \left[ (1+\lambda_m)^{-s/2} + \left(\frac{\eta}{n}\right)^{1/2} (1+\lambda_m)^{\frac{1}{2}(\tau-s)_+} \right]^2 + \sigma_\epsilon^2 n^{-1}. 
\end{align*} 
The hidden multiplicative constant is independent of $n$, $m$, $\eta$, $\delta$, and $\sigma_\epsilon^2$. 
\end{theorem}

Compared with the fixed-design result in Theorem~\ref{thm:bq_matern_fixed}, the random-design bound has a more concise form, at the cost of two additional conditions. 
We will show later that, for the Laplacian and Legendre operators, both conditions reduce to explicit constraints, and can be satisfied under specific tuning regimes. Unlike the fixed-design bound, the random-design bound does not involve the fill distance $h_n$. In fact, the quasi-uniform condition may not hold in random designs. For example, by Proposition 3.9 in \cite{reznikov:2015}, if ${\bfx_1,\ldots,\bfx_n}$ are i.i.d. from the uniform distribution on $\Omega$, then
$
    h_n = O_\PP\Big(\left(\frac{\log n}{n}\right)^{1/d}\Big),
$
which differs from the quasi-uniform rate $n^{-1/d}$ by a logarithmic factor.

\subsection{Extension to Gaussian-type weights}
\label{subsec:gaussian_theory}
We further extend the analysis to the Gaussian spectral weights defined in \eqref{eq:weights_gaussian}. In this case, we treat the lengthscale $\ell$ as a tunable parameter instead of a constant. Unlike the Mat\'ern case, the Gaussian kernel is not associated with a fixed finite Sobolev order, and the exponentially decaying spectral weights in \eqref{eq:weights_gaussian} require extra assumptions on the derivatives of eigenfunctions. The following result makes parameter interactions explicit for a fixed design.

\begin{theorem}
\label{thm:bq_gaussian_fixed} 
Suppose that Assumptions~\ref{assump:noise_vector},
\ref{assump:finite_density},
\ref{assump:f_smoothness},
\ref{assump:domain_padding},
\ref{assump:operator_coeff_smoothness_gaussian},
\ref{assump:local_elliptic},
\ref{assump:eigenpairs},
\ref{assump:eigenfunction_derivative}, and
\ref{assump:quasi_uniform} hold, and $w(\lambda_j)$ is the Gaussian weight function in \eqref{eq:weights_gaussian}.
Suppose $\ell, \eta>0$, and for the constant $C$ from Assumption~\ref{assump:eigenfunction_derivative},
\begin{align*}
     C^d (1+\lambda_m)^{d/2} \left[\log C + \log (1+\lambda_m) \right]^d
     \lesssim n
     \lesssim \frac{\eta}{\ell^d} (1+\lambda_m)^s,
\end{align*}
where the lower bound of $n$ holds with a sufficiently large multiplicative constant.
Then
\begin{align*} 
    \mathcal{Q}_n 
    &\lesssim h_n^{2s} + \left(1 + \lambda_m \right)^{-s}
    + \frac{\eta}{n\ell^d} 
    \max\left\{1, ~e^{\ell^2\lambda_m/2} (1+\lambda_m)^{-s} \right\}
    + \sigma_\epsilon^2 n^{-1} . 
\end{align*} 
The hidden multiplicative constant is independent of $n$, $m$, $\eta$, $\ell$, and $\sigma_\epsilon^2$. 
\end{theorem} 

As in the Mat\'ern case, the term
$
(1+\lambda_m)^{-s/2}
$
is the spectral truncation error. The Gaussian feature norm, however, introduces the exponential factor
$
\exp\left(\frac{1}{2}\ell^2\lambda_m\right).
$
This shows that increasing $m$ without simultaneously decreasing $\ell$ may amplify the finite-feature stability term, even though it reduces the truncation error. 
The theorem thus exposes a three-way trade-off: $m$ controls spectral truncation, $\ell$ controls both spectral growth and spatial resolution, and $\eta$ balances regularization bias against the magnitude of the quadrature weights.

For random designs, the same spectral structure remains, while deterministic fill-distance terms are replaced by concentration conditions. As defined in \eqref{eq:N.H}, $N(\Lambda)$ measures eigenvalue growth and $H(m)$ controls the pointwise concentration of eigenfunctions.

\begin{theorem}
\label{thm:bq_gaussian_random} 
Suppose that Assumptions~\ref{assump:noise_vector},
\ref{assump:bounded_density},
\ref{assump:f_smoothness},
\ref{assump:domain_padding},
\ref{assump:operator_coeff_smoothness_gaussian},
\ref{assump:local_elliptic},
\ref{assump:eigenpairs},
\ref{assump:eigenfunction_derivative}, and
\ref{assump:random_iid} hold, and $w(\lambda_j)$ is the Gaussian weight function in \eqref{eq:weights_gaussian}.
Suppose $\ell, \eta>0$, and $0<\delta<1$.
Let
$R_m=1+\lambda_m$ 
and 
$
L_p = \log \left( 2^{2p+4}N(2^{p+1}R_m) / \delta \right), ~p=0,1,\ldots.
$
For the constant $C$ from Assumption~\ref{assump:eigenfunction_derivative}, suppose that
\begin{align*}
    C^d R_m^{d/2}
    [\log(C R_m)]^d
    \log\left(\frac{4m}{\delta}\right)
    \lesssim
    n
    \lesssim
    \frac{\eta}{\ell^d}R_m^s,
\end{align*}
where the lower bound of $n$ holds with a sufficiently large implicit constant, and
\begin{align*}
    n
    \gtrsim
    \sum_{p=0}^\infty 2^{-sp}H\bigl(N(2^{p+1}R_m)\bigr)L_p.
\end{align*}
Then, with probability at least $1-\delta$ over $X_n$,
\begin{align*}
    \mathcal Q_n
    \lesssim 
    (1+\lambda_m)^{-s}
    +\frac{\eta}{n\ell^d}
       \max\left\{1,e^{\ell^2\lambda_m/2} (1+\lambda_m)^{-s}\right\}
    + \sigma_\epsilon^2 n^{-1} .
\end{align*}
The constants are independent of $n,m,\ell,\eta,\delta,\sigma_\epsilon^2$.
\end{theorem}

Although the error bound contains $\eta$ only through the increasing term
$
\eta/(n\ell^d),
$
decreasing $\eta$ does not necessarily improve the attainable rate. In fact, the sample-size condition
$
    n
    \lesssim
    \frac{\eta}{\ell^d}R_m^s
$
imposes an implicit lower bound on $\eta$.
Therefore, the choice of nugget $\eta$ reflects a genuine trade-off between reducing regularization error and preserving statistical stability.

\begin{remark}
    The bounds in Theorems~\ref{thm:bq_matern_random}
    and~\ref{thm:bq_gaussian_random} also hold when the design in Assumption~\ref{assump:random_iid} is instead sampled i.i.d. from a probability measure $\PP'$ on $\Omega$,
    which may differ from the integration measure $\PP$.
    Suppose that $\PP'$ admits a Lebesgue density $p'$ satisfying $ 0 < p_{\min}' \leq p'(\bfx) \leq p_{\max}' < \infty$, for all $\bfx \in \Omega$, for some constants $p_{\min}'$ and $p_{\max}'$, and that all remaining hypotheses are retained. Then the error bounds and sample-size conditions retain their form, with hidden multiplicative constants additionally depending on $p_{\min}'$ and $p_{\max}'$.
\end{remark}

\section{Uncertainty analysis}
\label{sec:posterior_variance}
We next study the uncertainty quantification for the finite-feature working Gaussian process model. From \eqref{eq:working_model_dist}, the GP posterior variance of the quadrature under the working model is
\begin{align*}
    \widehat\varrho_m^2
    &:=
    \Var\left[
    I(\bar f_m)\mid\Y_n
    \right] \\
    &=
    \frac{\sigma_\epsilon^2}{\eta} \left[\iint_{\Omega \times \Omega} \left[\widehat k_m(\bfx, \bfx^\prime) - \widehat \bfk_{n}^\top(\bfx) \left(\widehat \K_n+ \eta \I_n \right)^{-1} ~\widehat \bfk_n(\bfx^\prime) \right] \ud\PP(\bfx) \ud\PP(\bfx^\prime) \right].
\end{align*}

For every $\eta>0$, the posterior variance equals the regularized worst-case quadrature error over the finite-feature space. 
The key issue for frequentist coverage is whether the quadrature bias is negligible relative to the sampling uncertainty and whether the posterior variance correctly reflects that uncertainty. The following theorem formalizes these two requirements through standardized bias and variance-calibration conditions.

\begin{theorem} \label{thm:posterior_coverage}
    Suppose Assumptions~\ref{assump:noise_vector} and
    \ref{assump:f_smoothness} hold, with $\sigma_\epsilon^2>0$. 
    Let 
    \begin{align*}
        \bm\beta_n := \left(
        \widehat \K_n+\eta\I_n
        \right)^{-1} \widehat{\bm\mu}_n,
        \qquad
        \Delta_n(f) := I(f) - \bm\beta_n^\top \F_n, 
        \qquad
        \F_n = \left[f(\bfx_1), \ldots, f(\bfx_n)\right]^\top.
    \end{align*}
    For $t \in (0,1)$, let $I_t = \left[\widehat I_m(f) \pm z_{1-t/2} ~\widehat\varrho_m\right]$, where $z_{1-t/2}$ is the $(1-t/2)$-quantile of the standard normal distribution. Suppose $\bm\beta_n \neq \mathbf{0}$. If
    \begin{align} \label{eq:Delta.converge}
        \sup_{\|f\|_{H^s(\Omega)}\leq1} \frac{|\Delta_n(f)|}{\sigma_\epsilon \|\bm\beta_n \|_2} \rightarrow 0, ~\text{ as } n\rightarrow \infty,
    \end{align}
    then
    \begin{align*} 
        \liminf_{n\rightarrow \infty} \inf_{\|f\|_{H^s(\Omega)}\leq1} P_\epsilon\left(I(f) \in I_t \mid \X_n \right) \geq 1-t.
    \end{align*}
    Additionally, if 
    \begin{align} \label{eq:varrho.converge}
        \frac{\widehat\varrho_m}{\sigma_\epsilon  \|\bm\beta_n \|_2} \rightarrow 1,  ~\text{ as } 
 n\rightarrow \infty,
    \end{align}
    then
    \begin{align*}
        \sup_{\|f\|_{H^s(\Omega)}\leq1} 
        \left| P_\epsilon\left(I(f) \in I_t \mid \X_n \right) - (1-t) \right| \rightarrow 0,  ~\text{ as }  n \rightarrow \infty.
    \end{align*}
\end{theorem}

The quantity 
$\frac{|\Delta_n(f)|}{\sigma_\epsilon \|\bm\beta_n \|_2}$ in \eqref{eq:Delta.converge}
is the standardized bias for a given $f$. More precisely, the numerator $|\Delta_n(f)|$ quantifies the quadrature bias, while the denominator $\sigma_\epsilon\|\bm\beta_n\|_2$ is the conditional sampling standard deviation of the estimator. Thus, condition \eqref{eq:Delta.converge} means that the quadrature bias is uniformly negligible relative to the stochastic fluctuation caused by the observation noise. 
The first conclusion in Theorem~\ref{thm:posterior_coverage} guarantees that the worst-case coverage over the Sobolev unit ball does not asymptotically fall below the nominal level $1-t$, although the credible intervals may remain conservative. 

The quantity $\frac{\widehat\varrho_m}{\sigma_\epsilon  \|\bm\beta_n \|_2}$ in \eqref{eq:varrho.converge}
compares the posterior standard deviation with the true conditional sampling standard deviation and measures the frequentist calibration of the posterior uncertainty. The additional condition $\frac{\widehat\varrho_m}{\sigma_\epsilon  \|\bm\beta_n \|_2}\to1$ means that the posterior radius asymptotically matches the actual stochastic fluctuation of the estimator. Together with condition \eqref{eq:Delta.converge}, it yields the second conclusion, namely uniform asymptotically exact coverage. This second conclusion is stronger than the first one because it rules out both asymptotic undercoverage and asymptotic overcoverage. In the following section, we illustrate these two coverage results in Theorem~\ref{thm:posterior_coverage} and verify the corresponding bias condition \eqref{eq:Delta.converge} and variance-calibration condition \eqref{eq:varrho.converge} for specific combinations of spectral kernels and operator-adapted bases, including the Mat\'ern and Gaussian weights with concrete examples.

\section{Examples of eigenpairs and tuning rules}
\label{sec:concrete_bases}

The preceding theorems are stated for an abstract differential operator, and the error bounds depend on the eigensystem through quantities such as $\lambda_m$, $N(\Lambda)$, and $H(m)$. We now specialize the theory to two concrete operators, the Laplacian operator and the tensor-product Legendre operator, whose eigenpairs and basis integrals are analytically tractable. Such specializations convert the abstract error bounds and concentration conditions into explicit tuning rules for $m$, $\eta$, and, in the Gaussian case, $\ell$. 
In particular, under the feature-budget constraint $m \lesssim n$, we show that Mat\'ern weights attain the corresponding benchmark convergence rates for $\mathcal Q_n$ across fixed and random designs, noisy and noiseless observations, and both oversmoothing and undersmoothing regimes. 
For Gaussian weights, the resulting quadrature error attains the optimal rate in the noisy setting and the benchmark polynomial rate up to logarithmic factors in the noiseless setting.

Throughout this section, we take
$
L=B+1.
$
This choice keeps the target domain $\Omega$ uniformly separated from the artificial boundary of $D_L$, while treating $L$ as a constant independent of $n$. Consequently, factors depending only on $L$ do not affect the convergence rates.

\subsection{Laplacian eigenbasis}
The first choice is the Laplace operator
\begin{align} \label{eq:op_laplacian}
    \mathcal T_L = -\Delta = - \sum_{j=1}^d \partial_j^2.
\end{align}
The eigensystem of the Laplacian operator has different solutions under different boundary conditions:
\begin{align*}
\begin{array}{c|c|c|c}
\hline
\text{Boundary} 
& \text{Index } \bm\alpha 
& \lambda_{\bm\alpha} 
& \phi_{\bm\alpha}(\bfx) 
\\
\hline
\text{Dirichlet} 
& \NN_+^d 
& \displaystyle \left\|\frac{\pi \bm\alpha}{2 L}\right\|^2
& \displaystyle 
\prod_{p=1}^d 
\frac{1}{\sqrt{L}}
\sin\left(\frac{\alpha_p \pi(x_p+L)}{2L}\right)
\\[1.2em]
\text{Neumann} 
& \NN^d 
& \displaystyle \left\|\frac{\pi \bm\alpha}{2 L}\right\|^2
& \displaystyle 
\quad \prod_{p=1}^d 
c_{\alpha_p}
\cos\left(\frac{\alpha_p \pi(x_p+L)}{2L}\right),
~c_0=(2L)^{-1/2},
~c_r=L^{-1/2}, ~r \geq 1
\\[1.2em]
\text{Periodic} 
& \mathbb Z^d 
& \displaystyle \left\|\frac{\pi \bm\alpha}{ L}\right\|^2
& \displaystyle 
(2L)^{-d/2}
\exp\left(\frac{i\pi}{L} \bm\alpha^\top \bfx\right), ~ i=\sqrt{-1} 
\\[1.2em]
\hline
\end{array}
\end{align*}
The three boundary conditions produce different feature representations, but they have the same eigenvalue growth and uniform interior eigenfunction bounds. Therefore, the boundary condition affects the numerical form of the features and their integrals, but not the convergence orders derived below. For quadrature, the corresponding basis integrals reduce to trigonometric or exponential transforms of $\PP$, which are available for many commonly used probability measures. In particular, the basis associated with the Dirichlet boundary condition coincides with that used in the HSGP method in \cite{Solin:2020}, while the basis associated with the Periodic boundary coincides with that used in the EFGP method in \cite{Greengard:2025, Barnett:2024}. Compared with HSGP and EFGP, which construct multidimensional features through tensor-product or regular-grid structures, our finite-feature method selects basis functions directly according to the ordering of the operator eigenvalues. As a result, our method avoids an explicit tensor-product truncation across dimensions and can be more amenable to higher-dimensional settings.

\begin{remark}
    For complex-valued eigenfunctions such as the periodic Laplacian, let $\mathcal J_m\subset\mathbb Z^d$ denote the multi-index set of the retained eigenfunctions. 
    Then the truncated spectral kernel is understood as
    $
        \widehat{k}_m(\bfx, \bfx')
        =
        \sum_{\bm\alpha \in \mathcal J_m}
        w(\lambda_{\bm\alpha})
        \phi_{\bm\alpha}(\bfx) \overline{\phi_{\bm\alpha}(\bfx')}.
    $
    The index set can be chosen to be symmetric, namely
    $\bm\alpha \in \mathcal J_m$ implies $-\bm\alpha \in \mathcal J_m$, with $w(\lambda_{\bm\alpha})=w(\lambda_{-\bm\alpha})$. Pairing the frequencies $\bm\alpha$ and $-\bm\alpha$
    then gives a real-valued kernel. 
    Thus, the periodic case admits a real orthonormal eigenbasis, under which all statements in this paper remain unchanged.
\end{remark}

\subsection{Legendre eigenbasis}
The second choice is the classical tensor-product Legendre operator on $D_L$, defined by
\begin{align} \label{eq:op_legendre}
    \mathcal T_L 
    = 
    -\sum_{j=1}^d \partial_j \left[\left(1 - \frac{x_j^2}{L^2} \right) \partial_j \right].
\end{align} 
The self-adjoint realization and spectral properties of this operator are described in Proposition 1.6 and Theorem 3.6 of \cite{Littlejohn:2011}. The corresponding eigensystem $\{(\lambda_{\bm\alpha}, \phi_{\bm\alpha}) : \bm\alpha \in \NN^d \}$ is
\begin{align*}
    & \lambda_{\bm\alpha} = \frac{1}{L^2} \sum_{r=1}^d \alpha_r (\alpha_r + 1), \qquad
    \phi_{\bm\alpha}(\bfx) = \begin{cases}
        (2L)^{-d/2}, & \bm\alpha = \bm 0, \\
        \prod_{r=1}^d \sqrt{\frac{2\alpha_r + 1}{2L}} P_{\alpha_r}(x_r/L), & \bm\alpha \in \NN^d \setminus \{\bm 0\},
    \end{cases} \\
    & \text{where } P_r(z) := \frac{1}{2^r r!} \frac{d^r}{d z^r} (z^2 - 1)^r, \quad r = 0, 1, \ldots .
\end{align*}

This operator is uniformly elliptic on every interior domain containing $\Omega$, which is sufficient for the theory when $\Omega$ remains a fixed positive distance from $\partial D_L$. 
The Legendre basis is particularly useful for quadrature because its feature integrals reduce to polynomial moments of $\PP$, which are analytically available for many probability measures.

\subsection{Tuning rules for finite-feature quadrature}

We first apply Theorems~\ref{thm:bq_matern_fixed} and \ref{thm:bq_matern_random}. The resulting rates are identical for the Laplacian and Legendre bases because the general bounds depend only on their common eigenvalue and eigenfunction orders. The proposition also makes explicit how the admissible feature dimension $m$ and nugget parameter $\eta$ change between noisy and noiseless observations and between fixed and random designs.

\begin{proposition}
\label{prop:bq_matern_rates}
Let $\mathcal T_L$ be either the Laplace operator in \eqref{eq:op_laplacian}, or the tensor-product Legendre operator in \eqref{eq:op_legendre}. Let $L = B+1$.
Then Assumptions~\ref{assump:domain_padding}, \ref{assump:operator_coeff_smoothness}, \ref{assump:local_elliptic}, and \ref{assump:eigenpairs} hold, and the corresponding eigensystem satisfies
\[
    \lambda_m\asymp m^{2/d} L^{-2},
    \quad
    N(\Lambda)\lesssim L^d \Lambda^{d/2},
    \quad 
    H(m)\lesssim m L^{-d}, \quad \text{for all } \Lambda \geq 1,~ m \geq 2.
\]
Consider the Mat\'ern weight function $w(\lambda_j)$  in \eqref{eq:weights_matern}. 
Suppose Assumptions~\ref{assump:noise_vector}
and \ref{assump:f_smoothness} hold.
For the fixed-design cases, suppose
Assumptions~\ref{assump:finite_density}
and \ref{assump:quasi_uniform} hold.
For the random-design cases, suppose
Assumptions~\ref{assump:bounded_density}
and \ref{assump:random_iid} hold.
The tuning rule for different cases is as follows.
\begin{itemize}
\item Suppose noiseless observations and fixed design. If
\begin{align*}
    n^{\frac{1}{s}\min\{s,\tau\}}
    &\lesssim
    m
    \lesssim
    \begin{cases}
        n, & s<\tau,\\
        +\infty, & s\geq\tau,
    \end{cases}
    \\
    n^{-2\nu/d+2\min\{s,\tau\}/d}
    \left(
        n^{-2s/d}+m^{-2s/d}
    \right)
    &\lesssim
    \eta
    \lesssim
    n^{1-2\min\{s,\tau\}/d}
    m^{-2(\tau-s)_+/d},
\end{align*}
then
\[
    \mathcal Q_n
    \lesssim
    n^{-2\min\{s,\tau\}/d}.
\]

\item
Suppose noisy observations and fixed design. If
\begin{align*}
    n^{d/(2s)}
    \lesssim
    m
    \lesssim
    n^{\nu/(\tau-s)_+},
    \qquad
    n^{-2\nu/d}
    \lesssim
    \eta
    \lesssim
    m^{-2(\tau-s)_+/d},
\end{align*}
where the upper bound for $m$ is understood as $+\infty$ when
$\tau\leq s$, then
\[
    \mathcal Q_n
    \lesssim
    \left(1+\sigma_\epsilon^2\right)
    n^{-1}.
\]

\item
Suppose noiseless observations and random design. Let
$\delta\in(0,1)$ be fixed. 
If
\begin{align*}
    \eta
    \asymp
    n^{1-2\tau/d}
    \left[ \log (n/\delta) \right]^{2\tau/d},
    \qquad
    \left(
        \frac{n}{\log (n/\delta)}
    \right)^{\min\{s,\tau\}/s}
    \lesssim
    m
    \lesssim
    \frac{n}{\log (n/\delta)},
\end{align*}
then, with probability at least $1-\delta$ with respect to
$\X_n\overset{\mathrm{i.i.d.}}{\sim}\PP$,
\[
    \mathcal Q_n
    \lesssim
    \left(
        \frac{\log (n/\delta)}{n}
    \right)^{2\min\{s,\tau\}/d}.
\]

\item
Suppose noisy observations and random design. Let $\delta\in(0,1)$ be fixed. If
\begin{align*}
    n^{d/(2s)}
    &\lesssim
    m
    \lesssim
    \min\left\{
        \frac{n}{\log (n/\delta)},
        \,
        n^{\nu/(\tau-s)_+}
        \left[\log (n/\delta) \right]^{-\tau/(\tau-s)_+}
    \right\},
    \\
    n^{1-2\tau/d}
    \left[ \log (n/\delta) \right]^{2\tau/d}
    &\lesssim
    \eta
    \lesssim
    m^{-2(\tau-s)_+/d},
\end{align*}
where the second upper bound for $m$ is understood as $+\infty$
when $\tau\leq s$, then, with probability at least $1-\delta$,
\[
    \mathcal Q_n
    \lesssim
    \left(1+\sigma_\epsilon^2\right)
    n^{-1}.
\]
\end{itemize}
\end{proposition}

The noiseless rate depends on $\min\{s,\tau\}$, making the effect of smoothness mismatch explicit. 
When $s\leq \tau$, the resulting root mean squared error $n^{-s/d}$ agrees with the minimax deterministic quadrature rate over $H^s(\Omega)$, and with the rate established for exact Bayesian quadrature under quasi-uniform designs. When $s > \tau$, the rate $n^{-\tau/d}$ matches the standard RKHS benchmark of exact Bayesian quadrature with a $\tau$-smooth Mat\'ern kernel, as proved in \cite{Kanagawa:2020}. Compared with \cite{Kanagawa:2020}, our results additionally quantify the impact of the nugget term $\eta$, and include cases with random designs and noisy observations.
For noisy observations, both fixed and random designs attain the squared error rate $n^{-1}$. 
The root mean squared error rate $n^{-1/2}$ is faster than the rate $n^{-\frac{s}{2s+d}}$ established in full-kernel quadrature \cite{wynne:2021}.
Thus, a smaller feature dimension can retain the same accuracy while reducing the computational cost.
The random-design noiseless rate differs from its fixed-design counterpart only by logarithmic factors. These factors arise from concentration of the empirical feature system rather than from the finite-feature approximation itself.
For fixed quasi-uniform designs in the noiseless case, the sufficient range above forces
$m\asymp h_n^{-d}\asymp n$ when $s\leq\tau$. 
In all other cases covered by Proposition~\ref{prop:bq_matern_rates}, a choice
$m=o(n)$ is sufficient.
 
\begin{remark}
    Our tuning rule requires $L = B+1$ instead of $L \rightarrow \infty$. The reason is as follows. By Theorem~\ref{thm:bq_matern_fixed} and Proposition~\ref{prop:bq_matern_rates}, we have the
    term $(1 + \lambda_m)^{-s}$ in the upper bound of $\mathcal{Q}_n$ and $\lambda_m \asymp m^{2/d} L^{-2}$. If $L = L_n \rightarrow \infty$, we will suffer a loss of order $L^{2s}$. Unlike previous HSGP works \cite{Solin:2020} that require $L \rightarrow \infty$ to obtain accurate kernel approximation, we require $L = B + c_0$ for a constant $c_0 > 0$ to guarantee the stated quadrature convergence rates.
\end{remark}

We next specialize Theorems~\ref{thm:bq_gaussian_fixed} and \ref{thm:bq_gaussian_random} to the Laplacian and Legendre eigenbases. 
For both operators, the coefficient regularity condition in Assumption~\ref{assump:operator_coeff_smoothness_gaussian} holds for arbitrarily large $Q$, and their eigenfunctions are infinitely differentiable on $\Omega$. 
Moreover, Assumption~\ref{assump:eigenfunction_derivative} is directly verifiable from their explicit eigenpairs.
These properties allow the abstract Gaussian-weight bounds to be converted into explicit tuning rules for $m$, $\ell$, and $\eta$.

\begin{proposition}
\label{prop:bq_gaussian_rates}
Let $\mathcal T_L$ be either the Laplace operator in \eqref{eq:op_laplacian}, or the tensor-product Legendre operator in \eqref{eq:op_legendre}.
Let $L = B+1$.
Then Assumptions~\ref{assump:domain_padding}, \ref{assump:operator_coeff_smoothness_gaussian}, \ref{assump:local_elliptic}, and \ref{assump:eigenpairs} hold.
In addition, Assumption~\ref{assump:eigenfunction_derivative} holds for the two operators, and
\begin{align*}
    \lambda_m \asymp m^{2/d}L^{-2}, \quad 
    N(\Lambda)\lesssim L^d \Lambda^{d/2}, \quad H(m)\lesssim m L^{-d}, \quad \text{for all } \Lambda \geq 1,~ m \geq 2.
\end{align*}
Consider the Gaussian weight function $w(\lambda_j)$ in \eqref{eq:weights_gaussian}. 
Suppose Assumptions~\ref{assump:noise_vector}
and \ref{assump:f_smoothness} hold.
For the fixed-design cases, further suppose
Assumptions~\ref{assump:finite_density}
and \ref{assump:quasi_uniform} hold.
For the random-design cases, further suppose
Assumptions~\ref{assump:bounded_density}
and \ref{assump:random_iid} hold.
The tuning rule for different cases is as follows.
Constants in lower-bound constraints are taken sufficiently large, and constants in upper-bound constraints are sufficiently small, whenever required by Theorem~\ref{thm:bq_gaussian_fixed} and Theorem~\ref{thm:bq_gaussian_random}.
\begin{itemize}
\item
Suppose noiseless observations and fixed design. If
\[
    m \asymp \frac{n}{(\log n)^d},
    \qquad
    \ell \lesssim m^{-1/d} \sqrt{\log m},
    \qquad
    \eta \asymp n\ell^d m^{-2s/d},
\]
then
\[
    \mathcal Q_n
    \lesssim
    n^{-2s/d} (\log n)^{2s} .
\]

\item
Suppose noisy observations and fixed design. 
If
\begin{align*}
    n^{d/(2s)} \lesssim m \lesssim \frac{n}{(\log n)^d}, 
    \qquad
    \eta \lesssim m^{-1} (\log m)^{d/2}, 
\end{align*}
\begin{align*}
    \eta^{1/d} \lesssim \ell \lesssim \min\left\{m^{-1/d} \sqrt{\log m}, ~\left(\frac{\eta}{n} m^{2s/d} \right)^{1/d} \right\},
\end{align*}
then
\[
    \mathcal Q_n
    \lesssim
    \left(1+\sigma_\epsilon^2\right)
    n^{-1}.
\]

\item
Suppose noiseless observations and random design. Let
$\delta\in(0,1)$ be fixed. If
\begin{align*}
    m \asymp \frac{n}{\left[ \log (n/\delta) \right]^{d+1}},
    \qquad
    \ell \lesssim m^{-1/d} \sqrt{\log m},
    \qquad
    \eta \asymp n\ell^d m^{-2s/d} .
\end{align*}
Then, with probability at least $1-\delta$,
\[
    \mathcal Q_n
    \lesssim
    n^{-2s/d} ~\left[ \log (n/\delta) \right]^{2s+2s/d} .
\]

\item
Suppose noisy observations and random design. Let $\delta\in(0,1)$ be fixed. If
\begin{align*}
    n^{d/(2s)} &\lesssim m \lesssim \frac{n}{\left[ \log (n/\delta) \right]^{d+1}},
\end{align*}
\begin{align*}
    \eta \lesssim m^{-1} (\log m)^{d/2}, 
    \qquad
    \eta^{1/d} \lesssim \ell \lesssim \min\left\{m^{-1/d} \sqrt{\log m}, ~\left(\frac{\eta}{n} m^{2s/d} \right)^{1/d} \right\},
\end{align*}
then, with probability at least $1-\delta$,
\[
    \mathcal Q_n
    \lesssim
    (1+\sigma_\epsilon^2) n^{-1}.
\]
\end{itemize}
\end{proposition}

Under the Gaussian weight function, the error attains the optimal Sobolev polynomial rate, up to a logarithmic factor. This is expected, since the RKHS of Gaussian kernel contains infinitely smooth functions, while the unknown function is assumed to have only finite Sobolev regularity, $f\in H^s(\Omega)$ with $s>d/2$.

\begin{remark}
    The tuning rules in Propositions~\ref{prop:bq_matern_rates} and \ref{prop:bq_gaussian_rates} are oracle rules used to characterize the attainable statistical rates. In particular, some of them depend on the unknown target smoothness $s$. 
    If a lower smoothness bound $s_{\min}>d/2$ is known and the actual smoothness satisfies $s\ge s_{\min}$, then the conclusions of Propositions~\ref{prop:bq_matern_rates} and~\ref{prop:bq_gaussian_rates} can be applied with $s_{\min}$ in place of $s$. 
    In practice, one approach is to estimate $s$ using higher-order quadratic variation methods, available for both noiseless \citep{Loh:2021} and noisy \citep{Loh:2023} observations. Such estimates can be used to construct a plausible range for $s$. Alternatively, one may use a sufficiently smooth kernel with a large $\tau$ and adapt to the unknown function regularity through a data-dependent choice of the lengthscale $\ell$, for example via Bayesian updating procedures \citep{VV:2009, FangBhadra:2025}. 
    These approaches provide practical ways to relax the assumption that $s$ is known in advance.
\end{remark}

\subsection{Frequentist coverage of credible intervals}

We now apply Theorem~\ref{thm:posterior_coverage} to the Laplacian and Legendre eigenbases.
The preceding approximation bounds and tuning rules allow us to verify the standardized-bias condition for both the Mat\'ern and Gaussian spectral weights. For fixed quasi-uniform designs, this yields credible intervals whose asymptotic coverage is at least the nominal level. For random designs, concentration of the empirical feature Gram matrix further yields asymptotic variance matching and hence uniform asymptotically exact coverage.

\begin{corollary}[Fixed-design credible intervals]
\label{cor:fixed_design_coverage}
Suppose Assumptions~\ref{assump:noise_vector},
\ref{assump:finite_density},
\ref{assump:f_smoothness}, and
\ref{assump:quasi_uniform} hold, with
$\sigma_\epsilon^2>0$.
Let $I_t$ be defined as in
Theorem~\ref{thm:posterior_coverage}. 
Let $\mathcal T_L$ be either the Laplace operator in \eqref{eq:op_laplacian}, or the tensor-product Legendre operator in \eqref{eq:op_legendre}. 
Let $L=B+1$.

\begin{itemize}
    \item Suppose the Mat\'ern weights in \eqref{eq:weights_matern}
    are used, with $\tau=\nu+d/2$. Let
    \begin{align*}
        n^{d/(2s)} = o(m),
        \qquad
        m = o\left( n^{\nu/(\tau-s)_+} \right),
        \qquad
        n^{-2\nu/d}
        &\lesssim \eta
        =o\left(
            m^{-2(\tau-s)_+/d}
        \right),
    \end{align*}
    where $n^{\nu/(\tau-s)_+}$ is understood as $+\infty$ if $\tau \leq s$.

    \item Suppose the Gaussian weights in
    \eqref{eq:weights_gaussian} are used. Let
    \begin{align*}
        n^{d/(2s)} = o(m),
        \qquad
        m \lesssim \frac{n}{(\log n)^d},
    \end{align*}
    and
    \begin{align*}
        \eta^{1/d} = o(\ell),
        \qquad 
        \ell
        \lesssim
        \min\left\{
            m^{-1/d}\sqrt{\log m},
            ~\left(\frac{\eta}{n}m^{2s/d}\right)^{1/d}
        \right\}.
    \end{align*}
    In the Gaussian conditions, constants in upper-bound constraints are chosen sufficiently small whenever required by Theorem~\ref{thm:bq_gaussian_fixed}.
\end{itemize}
Under either set of conditions, 
\begin{align*}
    \sup_{\|f\|_{H^s(\Omega)}\leq1} \frac{|\Delta_n(f)|}{\sigma_\epsilon \|\bm\beta_n \|_2} \rightarrow 0, ~n\rightarrow \infty,
\end{align*}
Consequently, by Theorem~\ref{thm:posterior_coverage},
\begin{align*}
    \liminf_{n\rightarrow\infty}
    \inf_{\|f\|_{H^s(\Omega)}\leq1}
    P_\epsilon
    \left(
    I(f)\in I_t\mid\X_n
    \right)
    \geq1-t.
\end{align*}
\end{corollary}

\begin{corollary}[Random-design credible intervals]
\label{cor:random_design_coverage}
Suppose Assumptions~\ref{assump:noise_vector},
\ref{assump:bounded_density},
\ref{assump:f_smoothness}, and
\ref{assump:random_iid} hold, with
$\sigma_\epsilon^2>0$.
Let $I_t$ be defined as in
Theorem~\ref{thm:posterior_coverage}. 
Let $\mathcal T_L$ be either the Laplace operator in \eqref{eq:op_laplacian}, or the tensor-product Legendre operator in \eqref{eq:op_legendre}. 
Let $L=B+1$.
\begin{itemize}
    \item Suppose the Mat\'ern weights in \eqref{eq:weights_matern}
    are used, with $\tau=\nu+d/2$. Let
    \begin{align*}
        n^{d/(2s)} = o(m),
        \qquad
        m \lesssim \frac{n}{\log n},
        \qquad
        n^{1-2\tau/d} (\log n)^{2\tau/d}=o(\eta),
        \qquad
        \eta=o\left(m^{-2(\tau-s)_+/d} \right).
    \end{align*}

    \item Suppose the Gaussian weights in
    \eqref{eq:weights_gaussian} are used. Let
    \begin{align*}
        n^{d/(2s)} = o(m),
        \qquad
        m = o\left(\frac{n}{(\log n)^{d+1}}\right),
    \end{align*}
    and
    \begin{align*}
        \eta^{1/d} = o(\ell),
        \qquad 
        \ell
        \lesssim
        \min\left\{
            m^{-1/d}\sqrt{\log m},
            ~\left(\frac{\eta}{n}m^{2s/d}\right)^{1/d}
        \right\}.
    \end{align*}
    In the Gaussian conditions, constants in upper-bound constraints are chosen sufficiently small whenever required by Theorem~\ref{thm:bq_gaussian_random}.
\end{itemize}
Under either set of conditions, 
\begin{align*}
    n\|\bm\beta_n\|_2^2
    \xrightarrow{\PP} 1,
    \qquad
    \frac{n\widehat\varrho_m^2}{\sigma_\epsilon^2}
    \xrightarrow{\PP} 1,
\end{align*}
Consequently, by Theorem~\ref{thm:posterior_coverage},
\begin{align*}
    \sup_{\|f\|_{H^s(\Omega)}\leq1}
    \left|
    P_\epsilon
    \left(
    I(f)\in I_t\mid\X_n
    \right)
    -(1-t)
    \right|
    \xrightarrow{\PP} 0.
\end{align*}
\end{corollary}

The conditions in Corollaries~\ref{cor:fixed_design_coverage}
and \ref{cor:random_design_coverage} are closely related to the rate-optimal tuning rules in Propositions~\ref{prop:bq_matern_rates} and \ref{prop:bq_gaussian_rates}. 
In particular, for noisy observations, Propositions~\ref{prop:bq_matern_rates} and \ref{prop:bq_gaussian_rates} require the feature dimension to satisfy the lower bound $m\gtrsim n^{d/(2s)}$ to attain the benchmark squared-error rate $n^{-1}$. 
To verify the standardized-bias condition in Theorem~\ref{thm:posterior_coverage}, it is sufficient to strengthen this requirement only slightly, from $m\gtrsim n^{d/(2s)}$ to $n^{d/(2s)}=o(m)$.
Thus, asymptotically valid uncertainty quantification does not require a substantially larger feature budget than rate-optimal estimation, as it only requires $m$ to grow faster than the minimal sufficient order. The remaining conditions on $\eta$ and, for Gaussian weights, $\ell$ similarly strengthen the corresponding tuning ranges in Propositions~\ref{prop:bq_matern_rates} and \ref{prop:bq_gaussian_rates} so that the quadrature bias becomes negligible relative to the sampling standard deviation. 
For random designs, these conditions additionally ensure the variance-matching requirement in Theorem~\ref{thm:posterior_coverage}. 

\begin{remark}
The credible intervals in the two corollaries are oracle intervals because their posterior variances involve the true noise level $\sigma_\epsilon^2$. When $\sigma_\epsilon^2$ is unknown, it may be replaced by an estimator $\widehat\sigma_{\epsilon,n}^2$, such as the estimators in \cite{Liitiainen:2010} or \cite{Shen:2020Variance}. In particular, using prior scale $\widehat\sigma_{\epsilon,n}^2/\eta$ and working noise variance $\widehat\sigma_{\epsilon,n}^2$ preserves their ratio $\eta$ and leaves the posterior mean and the quadrature weights $\bm\beta_n$ unchanged. If $\widehat\sigma_{\epsilon,n}/\sigma_\epsilon\to1$ uniformly over the Sobolev unit ball, then the plug-in credible intervals satisfy the same asymptotic coverage conclusions, with the convergence interpreted conditionally for fixed designs and in $\PP$-probability for random designs. 
\end{remark}

\section{Feature-budget lower bounds}
\label{sec:feature_lower_bound}
Upper bounds alone do not determine whether the prescribed feature budget is necessary. 
For quadrature with a fixed uniform measure, no universal lower bound forcing the feature dimension to diverge can be derived solely from Sobolev smoothness and quasi-uniformity. In fact, a suitable design may encode an accurate quadrature rule using only a constant number of features. 
The situation changes, however, when the finite-feature procedure is required to attain a target rate uniformly over a nondegenerate class of integration measures. We show that, for fixed designs and densities uniformly bounded above and below, the worst-case quadrature error is bounded below by a constant multiple of $m^{-2s/d}$. Consequently, for the bounded density class, the lower bounds on $m$ in Propositions~\ref{prop:bq_matern_rates} and \ref{prop:bq_gaussian_rates} are not only sufficient but also necessary for attaining the corresponding optimal fixed-design quadrature rates.

\subsection{Feature-budget lower bound for non-degenerate density class}
\begin{proposition} \label{prop:density_class_lower_bound}
    Let $\Omega = (-B, B)^d$ and let $\X_n = \{\bfx_i\}_{i=1}^n \subset \Omega$ be a set of distinct points. 
    Suppose Assumptions~\ref{assump:noise_vector} and
    \ref{assump:f_smoothness} hold.
    Suppose $0 < p_{\min} < (2B)^{-d} < p_{\max} < \infty$.
    Let $\mathcal P = \{p(\bfx): \int_\Omega p(\bfx) \ud\bfx = 1, ~ p_{\min}\leq p(\bfx) \leq p_{\max}\}$.
    Then there exists a density $p(\bfx) \in \mathcal P$, depending on $m,n,\X_n$, such that
    \begin{align*}
        \sup_{\|f\|_{H^s(\Omega)}\leq 1} \EE_\epsilon\left[\left|I(f) - \widehat I_m(f) \right|^2 \middle| \X_n\right]
        \gtrsim
        m^{-2s/d}.
    \end{align*}
    The hidden multiplicative constant is independent of $n$, $m$, $\eta$, and $\sigma_\epsilon^2$. 
\end{proposition}

In Proposition~\ref{prop:bq_matern_rates}, under tuning rules for noiseless observations and fixed design, the lower bound of feature budget is
$m \gtrsim n^{\frac{1}{s} \min\{s,\tau\}}$. 
By Proposition~\ref{prop:density_class_lower_bound}, a necessary condition to obtain the optimal error rate is $n^{-2\min\{s,\tau\}/d} \gtrsim m^{-2s/d}$, thus $m \gtrsim n^{\frac{1}{s} \min\{s,\tau\}}$, so the lower bound $n^{\frac{1}{s} \min\{s,\tau\}}$ is also necessary. 
Similarly, the lower bounds on $m$ in other fixed design tuning rules in Propositions~\ref{prop:bq_matern_rates} and \ref{prop:bq_gaussian_rates} are also sufficient and necessary whenever a target error rate is required to hold uniformly over the density class $\mathcal P$.

\subsection{Why no point-wise feature-budget requirement holds for a fixed measure}
\begin{proposition} \label{prop:m1_design}
    Suppose Assumption~\ref{assump:noise_vector} holds.
    Suppose $\PP$ is the uniform probability measure on $\Omega = (-B, B)^d$. Let $m = 1$ and consider the rank-one constant working kernel
    \begin{align*}
        \widehat k_1(\bfx, \bfx') = a > 0.
    \end{align*}
    For any $s>d/2$, along an unbounded sequence of sample sizes $n$,
    there exist designs $\X_n$ satisfying
    $
        h_n \asymp q_n \asymp n^{-1/d},
    $
    with constants uniform over this sequence, such that, if
    $\eta\lesssim n^{1-s/d}$, then
    \begin{align*}
        \sup_{\|f\|_{H^s(\Omega)}\leq 1} \EE_\epsilon\left[\left|I(f) - \widehat I_1(f) \right|^2 \middle| \X_n\right]
        \lesssim
        n^{-2s/d} + \sigma_\epsilon^2 n^{-1}.
    \end{align*}
    The hidden multiplicative constant is independent of $n$, $m$, $\eta$, and $\sigma_\epsilon^2$. 
\end{proposition}

For Legendre or Laplacian-Neumann operator, when $m = 1$, the first eigenfunction is a constant $(2L)^{-d/2}$. Proposition~\ref{prop:m1_design} shows that there exists a specific quasi-uniform sequence such that the quadrature error achieves optimal rate with the rank-one working kernel.
Consequently, for the fixed integration measure considered there, no lower bound of the form $\mathcal Q_n \gtrsim m^{-2s/d}$ can hold uniformly over all quasi-uniform designs. When $m=1$, the right-hand side equals one, while $\mathcal Q_n$ converges to zero along the constructed sequence. This also illustrates that accurate quadrature can arise from a low-rank working kernel that need not accurately approximate the original kernel.

For a constant working kernel
$\widehat k_1(\bfx,\bfx')=a>0$ with fixed $a$, the estimator satisfies
$
    \widehat I_1(f)
    =
    \frac{an}{an+\eta}
    \frac1n\sum_{i=1}^n y(\bfx_i).
$
If
$\bfx_1,\ldots,\bfx_n\overset{\mathrm{i.i.d.}}{\sim}\PP$
and $\eta=o(n)$, then $\widehat I_1(f)\xrightarrow{\mathbb P}I(f)$ for each fixed
$f\in H^s(\Omega)$. Thus, one feature already suffices for convergence under random design. This result shows that for a fixed measure and random design, no universal lower bound forcing $m \to \infty$ can be derived.

\section{Consequences for related kernel estimators}
\label{sec:related_estimators}
The operator-adapted finite-feature construction is not specific to kernel quadrature. Since it provides a finite-dimensional representation of both the kernel matrix and the associated regularization norm, the same features can be used to construct computationally efficient estimators for other kernel-related problems. We first illustrate this extension through kernel ridge regression.

\subsection{Kernel ridge regression}
Define the finite-feature KRR estimator
\[
    \widehat f_m(\bfx)
    =
    \bm\phi(\bfx)^\top
    \left(
        \bm\Phi^\top \bm\Phi+n\lambda \W^{-1}
    \right)^{-1}
    \bm\Phi^\top\Y_n,
\]
where $\lambda = \eta/n >0$ is the regularization parameter. 
For the approximated kernel ridge estimator, the prediction error is defined as
\begin{align*} 
    \mathcal R_n
    :=
    \sup_{\|f\|_{H^s(\Omega)} \leq 1} \EE_\epsilon \left[ \left\|\widehat f_m - f \right\|_{L^2(\Omega)}^2 \mid \X_n\right].
\end{align*}
The following proposition shows that the operator-adapted finite-feature KRR estimator attains the Sobolev minimax prediction rate, matching the benchmark rate of optimally tuned kernel ridge regression, while using a reduced feature dimension.

\begin{proposition}[KRR rates for Mat\'ern-type weights \eqref{eq:weights_matern}]
\label{prop:krr_matern_rates}
Suppose Assumptions~\ref{assump:noise_vector}
and \ref{assump:f_smoothness} hold, with
$\sigma_\epsilon^2>0$.
Let $\mathcal T_L$ be either the Laplace operator in \eqref{eq:op_laplacian}
or the tensor-product Legendre operator in \eqref{eq:op_legendre},
with $L=B+1$.
For the fixed-design case, suppose
Assumption~\ref{assump:quasi_uniform} holds.
For the random-design case, suppose
Assumptions~\ref{assump:bounded_density}
and \ref{assump:random_iid} hold.
\begin{itemize}
\item
Fixed design. If
\[
    m\asymp n^{d/(2s+d)},
    \qquad
    n^{-2\tau/d}
    \lesssim
    \lambda
    \lesssim
    n^{-2\max\{s,\tau\}/(2s+d)},
\]
then
\[
    \mathcal R_n
    \lesssim
    \left(1+\sigma_\epsilon^2\right)
    n^{-2s/(2s+d)}.
\]

\item
Random design. For a fixed
$\delta\in(0,1)$, if
\[
    m\asymp n^{d/(2s+d)},
    \qquad
    \left(
        \frac{\log n}{n}
    \right)^{2\tau/d}
    \lesssim
    \lambda
    \lesssim
    n^{-2\max\{s,\tau\}/(2s+d)},
\]
then, with probability at least $1-\delta$ with respect to
$\X_n\overset{\mathrm{i.i.d.}}{\sim}\PP$,
\[
    \mathcal R_n
    \lesssim
    \left(1+\sigma_\epsilon^2\right)
    n^{-2s/(2s+d)}.
\]
\end{itemize}
\end{proposition}

Thus, in both the fixed- and random-design settings, the finite-feature KRR estimator attains the Sobolev minimax squared $L^2$-risk rate $n^{-2s/(2s+d)}$.
The conclusion holds in both the oversmoothing and undersmoothing regimes, with the admissible regularization range depending on $\max\{s,\tau\}$. In particular, the target is assumed only to belong to $H^s(\Omega)$ and need not lie in the Mat\'ern RKHS. 
For the Gaussian kernel, the analogous KRR result holds and attains the minimax Sobolev convergence rate.

\begin{remark}
    Our finite-feature method avoids saturation of full-kernel KRR when $s > 2\tau$. A well-known result in the KRR literature \cite{wang:2022, LiZhangLin:2023saturation} is that when $s > 2\tau$, the squared $L^2$-risk rate of full-kernel KRR with the Mat\'ern-$\nu$ kernel saturates at $n^{-\frac{4\tau}{4\tau+d}}$, instead of the minimax rate $n^{-\frac{2s}{2s+d}}$.
    The essential difference is that full-kernel KRR uses $\lambda$ both to reduce the shrinkage bias and to control the number of noisy directions. For a Mat\'ern kernel, \cite{LiZhangLin:2023saturation} establish lower bounds of order $\lambda^2$ for the squared bias and $n^{-1}\lambda^{-d/(2\tau)}$ for the variance. Balancing these lower bounds yields the saturation lower bound $n^{-4\tau/(4\tau+d)}$.
    Our finite-feature estimator separates these two roles. The cutoff $m$ controls the number of noisy directions, while $\lambda$ can be chosen sufficiently small that shrinkage within the retained feature space is negligible. The remaining squared bias is the omitted spectral tail of order $m^{-2s/d}$ and the variance is of order $m/n$. Balancing these terms gives $m\asymp n^{d/(2s+d)}$ and hence the minimax squared $L^2$-risk rate $n^{-2s/(2s+d)}$.
\end{remark}

\begin{proposition} \label{prop:function_recover_lower_bound}
    Let $\Omega=(-B,B)^d$ and $s>0$. Suppose $V_m$ is a fixed feature space and $\dimm(V_m) \leq m$. For any estimator $\tilde f_m \in V_m$,
    \begin{align*}
        \sup_{\|f\|_{H^s(\Omega)}\leq 1} \EE_\epsilon\left[\|f - \tilde f_m \|_{L^2(\Omega)}^2 \middle| \X_n \right] 
        \gtrsim 
        m^{-2s/d}.
    \end{align*}
\end{proposition}

The result is independent of the design,
noise distribution, regularization, and spectral weights, since these quantities do not enlarge the output space $V_m$. Consequently, attaining
a function recovery minimax rate $n^{-\frac{2s}{2s+d}}$ requires
$m\gtrsim n^{\frac{d}{2s+d}}$. This lower bound of feature budget $m$ matches the required order in Proposition~\ref{prop:krr_matern_rates}.

\subsection{Kernel mean embedding and maximum mean discrepancy}
The preceding quadrature results also imply convergence guarantees for kernel mean embeddings. Define
$
    \mu_{\PP}
    :=
    \int_\Omega
    k(\bfx,\cdot)\ud\PP(\bfx)
    \in\mathcal H_k,
$
and
$
    (w_1,\ldots,w_n)^\top
    =
    \Big(\widehat \K_n+\eta \I_n \Big)^{-1} \widehat{\bm\mu}_n.
$
\begin{proposition}[Kernel mean embedding]
\label{prop:kme_matern}
Let $k$ be the Mat\'ern kernel in \eqref{eq:kernel_matern}
with $\tau=\nu+d/2$.
Let $\mathcal T_L$ be either the Laplace operator in \eqref{eq:op_laplacian}
or the tensor-product Legendre operator in \eqref{eq:op_legendre},
with $L=B+1$.
For the fixed-design case, suppose
Assumptions~\ref{assump:finite_density}
and \ref{assump:quasi_uniform} hold.
For the random-design case, suppose
Assumptions~\ref{assump:bounded_density}
and \ref{assump:random_iid} hold.
Under the corresponding noiseless tuning
conditions of Proposition~\ref{prop:bq_matern_rates} with $s=\tau$,
\[
    \left\|\mu_{\mathbb P}
    -\sum_{i=1}^n w_i k(\mathbf x_i,\cdot)\right\|_{\mathcal H_k}^2
    \lesssim
    \begin{cases}
        n^{-2\tau/d}, & \text{fixed design},\\
        [n^{-1}\log(n/\delta)]^{2\tau/d}, & \text{random design}.
    \end{cases}
\]
The random-design bound holds with probability at least
$1-\delta$ over $\mathbf X_n$.
\end{proposition}

\section{Numerical experiments}
\label{sec:numerical_experiments}

We compare five quadrature methods on $\Omega=(-B,B)^d$, $d=2$.
The five methods are: full-kernel quadrature \eqref{eq:bq_predictor}, Laplacian finite-feature quadrature, Legendre finite-feature quadrature \eqref{eq:finite_feature_estimator}, Nystr\"om sampling for kernel integral \cite{Chatalic:2022, Chatalic:2025}, and Nystr\"om finite-feature quadrature for kernel matrix \cite{Rudi:2015, Rudi:2017FALKON}. 
Throughout, $B=1$, $L = B+1 = 2$.
The test function is the additive cusp-type Sobolev function
\begin{align*}
    f(\bfx)
    =
    \sum_{j=1}^d \frac{1}{j} \left|x_j - \frac{B}{2} \sin j \right|^q,
\end{align*}
which has a closed-form integral under both uniform and beta distributions on $\Omega$.
Since $f \in H^{q+\epsilon}(\Omega), ~\epsilon \in [0, 0.5)$ and $f \notin H^{q+0.5}(\Omega)$ when $q \notin 2\NN$,
we set $s = q+1/2-0.001$. 
For noisy experiments, the observations are
\[
    y_i=f(\bfx_i)+\varepsilon_i,
    \qquad
    \varepsilon_i\sim N(0,10^{-2}),
    \qquad
    i = 1, \ldots, n.
\]
The sample sizes are
$
    n\in\{64, 100, 144, 196, 256, 324, 400, 625, 841, 1024\}.
$
All random experiments are based on 100 independent replications, while deterministic estimators under fixed, noiseless designs are evaluated once.
For the fixed-design experiments, the design points are generated by tensor grids, and the integration measure is the uniform probability measure on $\Omega$.
For the random-design experiments, the observation points are sampled independently from the integration measure, whose density on $\Omega$ is
\begin{align*} 
    p_c(\bfx) &= (1-c) \prod_{j=1}^d \rho_{a,b}(x_j) + c \frac{1}{(2B)^d}, \\
    \rho_{a,b}(t)
    &=
    \frac{(t+B)^{a-1}(B-t)^{b-1}}
    {(2B)^{a+b-1}\mathrm{B}(a,b)}.
\end{align*}
Here $\rho_{a,b}$ is the density of a Beta$(a,b)$ random variable rescaled to $(-B,B)$, and $\mathrm{B}(a,b)$ denotes the Beta function.
We set $a=b=2$ and $c=0.1$. 

We first consider the Mat\'ern kernel \eqref{eq:kernel_matern} with smoothness parameter $\nu=2.5$, lengthscale $\ell=0.5$, and variance $\sigma^2=1$.
Writing $\tau=\nu+d/2$, we consider $q=\tau$ and $q=\tau-0.75$, corresponding to kernel undersmoothing ($s>\tau$) and oversmoothing ($s<\tau$), respectively, where $s=q+1/2-0.001$.
For Laplacian and Legendre BQ, we choose parameters based on the tuning rates in Proposition~\ref{prop:bq_matern_rates}, that is, 
\begin{align*}
    \eta = C_\eta
    \begin{cases}
        n^{-2\nu/d}, &\text{fixed design}, \\
        n^{-2\nu/d}(\log n)^{2\tau/d}, &\text{random design},
    \end{cases}
\end{align*}
\begin{align*}
    m = 
    \begin{cases}
        \left\lceil C_m n^{\min\{s,\tau\}/s} \right\rceil, &\text{fixed design, noiseless}, \\
        \left\lceil C_m \left(n/\log n\right)^{\min\{s,\tau\}/s}\right\rceil , &\text{random design, noiseless}, \\
        \left\lceil C_m \sqrt{n}\right\rceil , &\text{noisy}.
    \end{cases}
\end{align*}

Full-kernel BQ uses the same nugget $\eta$ as the finite-feature methods.
Nystr\"om-matrix BQ uses the same feature budget $m$ and nugget $\eta$, with $m$ landmarks sampled with replacement from the observation points with probabilities proportional to the ridge leverage scores
$
    \tau_i^{(n)}(\lambda)
    =
    \left[
        \K_n(\K_n+n\lambda\I_n)^{-1}
    \right]_{ii},
    ~\lambda=\eta/n.
$
The Mat\'ern kernel means required by full-kernel BQ and Nystr\"om-matrix BQ are approximated using randomized quasi-Monte Carlo with $2^{12}$ scrambled Sobol points and replicate-specific scrambling \cite{Sobol:1967}.

Nystr\"om-integral forms an unlabeled pool
$\Z_N=\{\bfz_i^{\mathrm{pool}}\}_{i=1}^N$, with
$N=\min\{10n,10^4\}$, sampled independently from the target integration measure.
This method selects $n$ evaluation nodes with replacement from this pool, with ridge score
$
    \tau_i^{(N)}(\lambda)
    =
    \left[
        \K_N(\K_N+N\lambda\I_N)^{-1}
    \right]_{ii},
$
$\lambda=19 \frac{\log(640N)}{N}$.
The estimator uses the empirical kernel mean
$
    \widehat\mu_N(\bfz)
    =
    \frac{1}{N}\sum_{i=1}^N k(\bfz, \bfz_i^{\mathrm{pool}})
$
and the Moore--Penrose pseudoinverse of the kernel matrix at the selected nodes, without an added estimator nugget.

\begin{figure}[t]
    \centering
    \includegraphics[width=\linewidth]{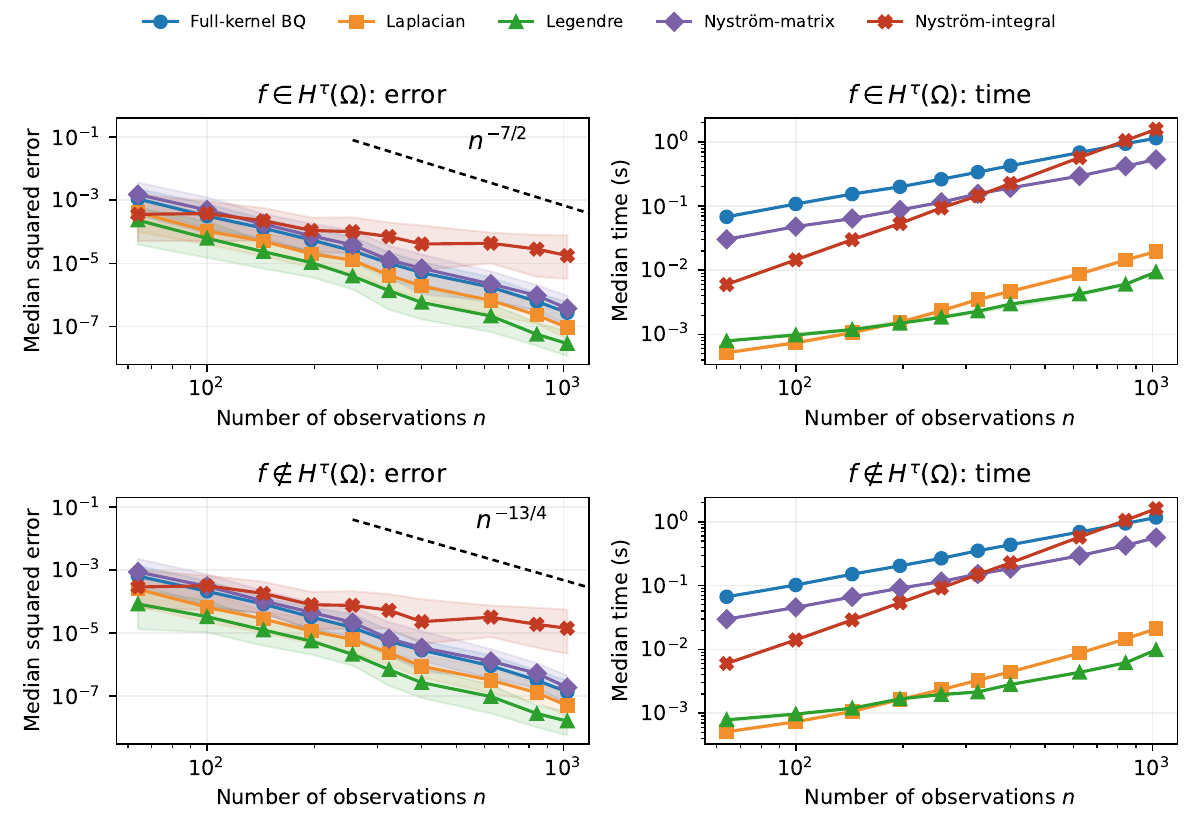}
    \caption{
        Random-design results for the noiseless Mat\'ern experiments in
        dimension $d=2$. The top and bottom rows correspond to
        $f\in H^\tau(\Omega)$ and $f\notin H^\tau(\Omega)$,
        respectively. 
        The left panels report the empirical median squared integration error, while the right panels report the median runtime.
        The black dashed lines show the limiting theoretical rates $n^{-2\min\{q+1/2, \tau\}/d}$ from Proposition~\ref{prop:bq_matern_rates}.
        The shaded bands show the 25th–75th percentiles.
    }
    \label{fig:matern-random-noiseless}
\end{figure}

Figure~\ref{fig:matern-random-noiseless} shows that the
Legendre method achieves the smallest error, followed by the Laplacian method, while both finite-feature methods are substantially faster than the kernel-based alternatives. 
Nystr\"om-matrix has errors close to full-kernel BQ and offers a moderate runtime advantage. 
Nystr\"om-integral shows slower error reduction and loses its initial runtime advantage as $n$ increases.
The method ranking remains similar inside and outside the kernel RKHS, with no apparent deterioration caused by leaving the kernel RKHS. 
The reference lines represent the predicted polynomial rates $n^{-2\tau/d}$ and $n^{-(2q+1)/d}$, respectively, omitting logarithmic factors. 
The steep error reduction of the finite-feature methods is qualitatively consistent with these benchmarks.

The finite-feature methods benefit from closed-form feature means and a computational cost of $O(nm^2+m^3)$. The errors from Laplacian and Legendre operators have the same orders, and these lower errors may also reflect beneficial spectral truncation and approximation constants suited to these test functions. 
Nystr\"om-integral's weaker performance may partly reflect its empirical kernel-mean approximation and its use of a pseudoinverse without a tuned estimator nugget. The ridge parameter used for node sampling does not regularize that estimator. 
Table~\ref{tab:bq-d2-relative} summarizes the remaining Mat\'ern-kernel experiments.

\begin{table}[t]
\centering
\footnotesize
\setlength{\tabcolsep}{5pt}
\renewcommand{\arraystretch}{1.13}
\caption{Relative Mat\'ern BQ results at $n=1024$. Each entry reports error and time ratios. The error or time ratio is the method's median squared error or computation time divided by that of full-kernel BQ. A smaller ratio is better.}
\label{tab:bq-d2-relative}
\begin{tabular}{@{}lcccc@{}}
\toprule
Setting & Laplacian & Legendre & \shortstack{Nystr\"om-matrix} & \shortstack{Nystr\"om-integral} \\
\midrule
Fixed, noiseless, $f\in H^\tau$
& $9.46\times10^{-3}\,/\,0.0328$
& $3.25\times10^{-4}\,/\,0.0115$
& $3.88\,/\,0.602$
& $2.52\times10^{3}\,/\,2.46$ \\

Fixed, noiseless, $f\notin H^\tau$
& $7.31\times10^{-3}\,/\,0.0288$
& $9.53\times10^{-3}\,/\,0.0102$
& $3.86\,/\,0.618$
& $3.84\times10^{3}\,/\,2.45$ \\

Fixed, noisy, $f\in H^\tau$
& $0.986\,/\,0.0317$
& $0.996\,/\,0.0129$
& $0.902\,/\,0.633$
& $37.7\,/\,2.45$ \\

Fixed, noisy, $f\notin H^\tau$
& $0.968\,/\,0.0317$
& $0.978\,/\,0.0130$
& $0.885\,/\,0.634$
& $36.0\,/\,2.45$ \\

Random, noiseless, $f\in H^\tau$
& $0.338\,/\,0.0172$
& $0.105\,/\,8.22\times10^{-3}$
& $1.36\,/\,0.465$
& $64.0\,/\,1.38$ \\

Random, noiseless, $f\notin H^\tau$
& $0.371\,/\,0.0180$
& $0.119\,/\,8.41\times10^{-3}$
& $1.40\,/\,0.485$
& $105\,/\,1.37$ \\

Random, noisy, $f\in H^\tau$
& $0.803\,/\,0.0184$
& $0.770\,/\,8.47\times10^{-3}$
& $1.00\,/\,0.497$
& $16.5\,/\,1.37$ \\

Random, noisy, $f\notin H^\tau$
& $0.884\,/\,0.0184$
& $0.782\,/\,8.49\times10^{-3}$
& $0.942\,/\,0.497$
& $16.5\,/\,1.37$ \\
\bottomrule
\end{tabular}
\end{table}

\section{Discussion}
\label{sec:discussion}
This work develops an operator-adapted finite-feature framework for scalable kernel quadrature and related kernel estimators. The main contribution is an end-to-end analysis that treats feature integration and finite-dimensional regularized inference within the same operator-adapted construction. The construction uses the eigensystem of a second-order differential operator on an enlarged domain. 
We establish convergence rates and frequentist coverage guarantees for operator-adapted finite-feature quadrature under fixed and random designs, allowing noisy observations and kernel smoothness misspecification.
We apply the resulting theory to both Laplacian and Legendre-type bases and cover fixed and random designs, noisy and noiseless observations, Mat\'ern and Gaussian type spectral weights, and matched and misspecified Sobolev smoothness regimes. 

Our results clarify when computational savings can be obtained without sacrificing statistical accuracy. The finite-feature quadrature estimators require solving an $m$-dimensional linear system with computational cost $O(nm^2+m^3)$, reducing the $O(n^3)$ cost of inverting an $n\times n$ kernel matrix. Under the tuning regimes established in this paper, one can take $m\lesssim n$, while retaining the optimal rate in the noisy regimes and the optimal or near-optimal rate in the noiseless regimes. The same operator-adapted construction also yields minimax-optimal finite-feature KRR estimators and corresponding guarantees for kernel mean embedding. 

Several extensions remain open. The present analysis assumes a bounded domain separated from the boundary of the auxiliary operator domain and relies on second-order operators with sufficiently regular coefficients. Extending the theory to irregular domains, nonstationary kernels, and operators adapted to more general geometries would broaden the range of analytically tractable feature systems. In addition, the parameter-tuning rules depend on the target Sobolev smoothness $s$. 
Developing data-adaptive smoothness and feature-selection rules and extending the analysis to more general spectral bases and integration measures are important directions for future work.

\newpage

\bibliographystyle{plainnat}
\bibliography{bibliography}

\end{document}